\documentclass[11pt]{amsart}

\usepackage{geometry}
\usepackage{todonotes}
\usepackage{amsfonts, amssymb, amscd}
\numberwithin{equation}{section}

\usepackage[symbol]{footmisc}

\usepackage{bm}
\usepackage{verbatim}
\usepackage{mathrsfs}
\usepackage{graphicx}
\usepackage{tikz-cd}
\usepackage{subcaption}
\usepackage{listings}
\usepackage{subfiles}
\usepackage[toc,page]{appendix}
\usepackage{mathtools}
\usepackage{comment}
\usepackage{enumerate}
\usepackage{enumitem}
\usepackage[all]{xy}

\usepackage{graphicx}
\graphicspath{{images/}}

\usepackage{appendix}
\usepackage{hyperref}
\hypersetup{
    colorlinks=true,
    citecolor=red,
    linkcolor=blue,
    filecolor=magenta,      
    urlcolor=red,
}
\newcommand{\bb}{\bm{b}}
\newcommand{\Mm}{{\bf{M}}}
\newcommand{\Bb}{{\bf{B}}}

\newcommand{\Dd}{{\bf{D}}}

\newcommand{\Cc}{\mathbb{C}}

\newcommand{\Qq}{\mathbb{Q}}

\newcommand{\Rr}{\mathbb{R}}

\newcommand{\Exc}{\operatorname{Exc}}

\newcommand{\rk}{\operatorname{rank}}

\newcommand{\loc}{\operatorname{loc }}

\newcommand{\Weil}{\operatorname{Weil}}

\newcommand{\Supp}{\operatorname{Supp}}

\newcommand{\mult}{\operatorname{mult}}

\newcommand{\Aa}{{\bf{A}}}

\newcommand{\Ff}{\mathcal{F}}
\newcommand{\Gg}{\mathcal{G}}

\newcommand{\Ii}{\Gamma}

\newcommand{\Ee}{\mathcal{E}}

\newcommand{\Sing}{\mathrm{Sing}}

\newcommand{\rank}{\mathrm{rank}}

\newtheorem{thm}{Theorem}[section]
\newtheorem{conj}[thm]{Conjecture}

\newtheorem{lem}[thm]{Lemma}
\newtheorem{prop}[thm]{Proposition}

\theoremstyle{definition}
\newtheorem{defn}[thm]{Definition}

\theoremstyle{definition}
\newtheorem{rem}[thm]{Remark}

\newtheorem{deflem}[thm]{Definition-Lemma}
\newtheorem{defthm}[thm]{Definition-Theorem}
\newtheorem{setup}[thm]{Set-up}

\newtheorem{nota}[thm]{Notation}

\theoremstyle{definition}

\begin{document}

\title{On global ACC for foliated threefolds}
\author{Jihao Liu, Yujie Luo, and Fanjun Meng}

\date{\today}

\begin{abstract}
In this paper, we prove the rational coefficient case of the global ACC for foliated threefolds. Specifically, we consider any lc foliated log Calabi-Yau triple $(X,\Ff,B)$ of dimension $3$ whose coefficients belong to a set $\Ii$ of rational numbers satisfying the descending chain condition, and prove that the coefficients of $B$ belong to a finite set depending only on $\Ii$.

To prove our main result, we introduce the concept of generalized foliated quadruples, which is a mixture of foliated triples and Birkar-Zhang's generalized pairs. With this concept, we establish a canonical bundle formula for foliations in any dimension.

As for applications, we extend Shokurov's global index conjecture in the classical MMP to foliated triples and prove this conjecture for threefolds with nonzero boundaries and for surfaces. Additionally, we introduce the theory of rational polytopes for functional divisors on foliations and prove some miscellaneous results.
\end{abstract}

\address{Department of Mathematics, Northwestern University, 2033 Sheridan Road, Evanston, IL 60208, USA}
\email{jliu@northwestern.edu}

\address{Department of Mathematics, Johns Hopkins University, 3400 N. Charles Street, Baltimore, MD 21218, USA}
\email{yluo32@jhu.edu}

\address{Department of Mathematics, Johns Hopkins University, 3400 N. Charles Street, Baltimore, MD 21218, USA}
\email{fmeng3@jhu.edu}

\maketitle

\tableofcontents

\section{Introduction}\label{sec:Introduction}
We work over the field of complex numbers $\mathbb C$.

The theory of \emph{foliations} is an intriguing subject that connects various areas of mathematics, including algebraic geometry, complex analysis, and topology. Foliations exhibit favorable properties in birational geometry and play a critical role in the minimal model program (Mori's program), particularly in Miyaoka's proof of the abundance conjecture in dimension three \cite{Miy87} (also see \cite[Chapter 9]{Kol+92}). In recent years, it has become apparent that many results in classical birational geometry extend to foliations, particularly in low dimensions. Specifically, the foundations of the minimal model program for foliated surfaces (cf. \cite{McQ08,Bru15}) and foliated threefolds (cf. \cite{CS20,Spi20,CS21,SS22}) have been established.

However, not all properties of conventional varieties hold for foliations. For example, the abundance conjecture fails in general, even for surfaces \cite[Theorem 3 IV.5.11]{McQ08}, and effective birationality also fails in general \cite[Paragraph before Theorem 1.4]{SS23}. As a result, the boundedness of foliations of general type can only be proved for surfaces under certain additional assumptions \cite{Che21a,Che21b,HL21,SS23}.

Despite these differences, we observe that the properties that do not hold for foliations are all ``positivity" properties, such as abundance, effective birationality, and boundedness of foliations of general type. In contrast, we still expect the ``negativity" or ``non-positivity" properties (e.g., Fano type, log Calabi-Yau) for foliations to behave similarly to conventional varieties. An important piece of evidence for this is the famous result of Campana-P\u{a}un \cite[Theorem 1.2]{CP19}, which suggests that any smooth projective variety has a non-pseudo-effective canonical class as long as it has a foliation structure with a non-pseudo-effective canonical class. This result has been recently generalized to the case of klt pairs \cite[Theorem 1.2]{ACSS21}.

With these ideas in mind, we prove some results on the ``negativity" or ``non-positivity" properties of foliations in this paper.

\medskip

\noindent\textbf{Global ACC for foliated threefolds}. The first main result of this paper is the global ACC and the index theorems for foliated log Calabi-Yau pairs in dimension $\leq 3$.

\begin{thm}\label{thm: global acc threefold}
Let $\Ii\subset [0,1]\cap\mathbb Q$ be a DCC set. Then there exists a positive integer $I$ depending only on $\Ii$ satisfying the following. Assume that $(X,\Ff,B)$ is a projective lc foliated triple such that $\dim X\leq 3$, $K_{\Ff}+B\equiv 0$, and $B\in\Ii$. Then:
\begin{enumerate}
    \item The coefficients of $B$ belong to the finite set $\frac{1}{I}\mathbb N\cap [0,1]$.
    \item If $B\not=0$, then $I(K_{\Ff}+B)\sim 0$.
\end{enumerate}
\end{thm}
We now summarize some previously established results related to Theorem \ref{thm: global acc threefold}.

\begin{enumerate}
\item When $\dim X = 2$, Theorem \ref{thm: global acc threefold}(1) has been proven in \cite[Theorem 2.5]{Che22}.
\item Under the hypotheses of Theorem \ref{thm: global acc threefold}, \cite[Theorem 10.1]{CS20} and \cite[Theorem 12.1]{CS21} prove that $K_{\Ff} + B \sim_{\mathbb Q} 0$ when $B$ is a $\Qq$-divisor.
\item When $\dim X = 2$ and $B=0$, Theorem \ref{thm: global acc threefold}(2) has been proven in \cite[Theorem 1]{Per05} when $\Ff$ is canonical, and in full generality in Theorem \ref{thm: surface index theorem global intro} below.
\item When $\Ff=T_X$, Theorem \ref{thm: global acc threefold}(1) is a special case of the well-known global ACC for usual pairs, which is proven in \cite[Theorem 1.5]{HMX14}. Theorem \ref{thm: global acc threefold}(2) was recently established in dimension $3$ in \cite[Theorem 13]{Xu19} as an application of \cite[Corollary 1.7]{Jia21}.
\end{enumerate}

We expect Theorem \ref{thm: global acc threefold}(2) to hold even when $\dim X=3$ and $B=0$, although the strategies to prove it are expected to be quite different. We would like to point out that Theorem \ref{thm: global acc threefold}(2) is not known even in the case when $\Ff=T_X$ and $B=0$ in dimension $\geq 4$. This highlights the difficulty of the problem and motivates further research in this direction.

\medskip

\noindent\textit{Idea of the proof of Theorem \ref{thm: global acc threefold}}. The proof of Theorem \ref{thm: global acc threefold} takes a different approach from the proof of the global ACC for usual pairs \cite[Theorem 1.5]{HMX14}. This is because the proof of \cite[Theorem 1.5]{HMX14} relies on the effective birationality for pairs of log general type \cite[Theorem 1.3(3)]{HMX14}, which fails in general for foliations \cite[Paragraph before Theorem 1.4]{SS23}. Additionally, since most foliations are not klt in codimension $1$, the foliation version of \cite[Theorem B]{HMX14} is also not useful when proving the global ACC. We remark that O. Das informed us that he and his collaborators have got some results related to the foliated version of \cite[Theorem B]{HMX14}. 

Therefore, an alternative approach is needed. A key observation is that any non-pseudo-effective foliation always has a non-trivial sub-foliation that is algebraically integrable (cf. Theorem \ref{thm: subfoliation algebraic integrable}). This observation indicates that non-trivial fibration structures are easier to obtain compared to the case of usual pairs. With a fibration structure, the goal is to show that a fixed multiple of $K_{\Ff}+B$ is linearly equivalent to the pullback of a Cartier divisor from the base of the fibration. To achieve this goal, there are three parts of the work to do:

\begin{enumerate}
\item Establish a canonical bundle formula for foliated fibrations so that we can control the structure of the base of the fibration. The details of this step are given in Theorem \ref{thm: cbf ftriple}.
\item When $B$ is a $\Qq$-divisor, show that the foliated log canonical divisor on the base has bounded index.
\item Reduce to the case when $B$ is a $\Qq$-divisor.
\end{enumerate}
We will discuss these three parts in the following sections.

\medskip

\noindent\textbf{Canonical bundle formulas and generalized foliated quadruples}. The first key ingredient in the proof of Theorem \ref{thm: global acc threefold} is the canonical bundle formula for foliations. While some special cases of the canonical bundle formula for foliations were introduced and used in \cite[Lemma 9.1]{CS21} and \cite[Proposition 3.6]{ACSS21}, we establish a more general canonical bundle formula for foliations by introducing the notion of \emph{generalized foliated quadruples}. This concept is a mixture of foliated triples and Birkar-Zhang's generalized pairs \cite{BZ16}, which allows us to extend and generalize the canonical bundle formula to foliations of any dimension. The proof of the canonical bundle formula for foliated triples involves several technical steps, including a careful analysis of the singularities of generalized foliated quadruples, but ultimately leads to a powerful tool for studying foliations and related fibration structures.

\begin{defn}
A \emph{generalized foliated quadruple (gfq)} $(X,\Ff,B,\Mm)/U$ consists of a normal quasi-projective variety $X$, a foliation $\Ff$ on $X$, an $\mathbb R$-divisor $B\geq 0$ on $X$, a projective morphism $X\rightarrow U$, and a $\bb$-divisor $\Mm$ on $X$ satisfying the following conditions: $\Mm=\sum m_i\Mm_i$ where $m_i\geq 0$ and $\Mm_i$ are nef$/U$ $\bb$-Cartier $\bb$-divisors, and $K_{\Ff}+B+\Mm_X$ is $\mathbb R$-Cartier.
\end{defn}

The singularities of gfqs are defined similarly to the singularities of usual pairs, generalized pairs, and foliated triples (cf. Definition \ref{defn: gfq singularity}). We establish the following canonical bundle formula for foliated triples.

\begin{thm}\label{thm: cbf ftriple}
Let $(X,\Ff,B)/U$ be an lc foliated triple and let $f: (X,\Ff,B)\rightarrow Z$ be an lc-trivial fibration$/U$ (cf. Definition \ref{defn: lc trivial fibration foliation}). Then there exists a gfq induced by a canonical bundle formula$/U$ (cf. Definition-Lemma \ref{deflem: discriminant and moduli part}) for $f: (X,\Ff,B)\rightarrow Z$. Moreover, for any gfq $(Z,\Ff_Z,B_Z,\Mm^Z)/U$ induced by a canonical bundle formula$/U$ for $f: (X,\Ff,B)\rightarrow Z$, we have:
\begin{enumerate}
\item $(Z,\Ff_Z,B_Z,\Mm^Z)$ is lc,
\item for any component $D$ of $B_Z$, 
$$\mult_DB_Z=1-\sup\{t\mid (X,\Ff,B+tf^*D)\text{ is lc over the generic point of }D\},\text{ and}$$
\item the class of $\Mm^Z$, up to $\Rr$-linear equivalence on high models, only depends on $(X,B)$ over the generic point of $Z$. In particular, for any integer $n$ such that $n(K_X+B)\sim 0$ over the generic point of $Z$, we may choose $\Mm^Z$ such that $n(K_{\Ff}+B)\sim nf^*(K_{\Ff_Z}+B_Z+\Mm^Z_Z)$.
\end{enumerate}
\end{thm}
We remark that Theorem \ref{thm: cbf ftriple} holds for all dimensions. For a detailed version of Theorem \ref{thm: cbf ftriple} and additional information on the case when $\dim X \leq 3$, we refer to Proposition \ref{prop: low-dimension cbf}.

We expect that the minimal model program for foliations also works for gfqs, especially in dimension $\leq 3$, such as cone theorems, contraction theorems, existence of flips, termination of flips, and the ACC for lc thresholds. We also expect Theorems \ref{thm: global acc threefold} and \ref{thm: cbf ftriple} to hold for gfqs in full generality. Due to technical reasons, we will not include these results in this paper, but we plan to explore them in future work. 

It is crucial to emphasize that the concept of gfqs plays a vital role in proving Theorem \ref{thm: global acc threefold}. The main reason behind this lies in the fact that even if $\Mm^Z$ is semi-ample (which is conjecturally true for gfqs induced by canonical bundle formulas and is known in low dimensions; see Proposition \ref{prop: low-dimension cbf}), it is generally not possible to find $0\leq G\sim_{\mathbb R}\Mm^Z_Z$ such that $(Z,\Ff_Z,B_Z+G)$ is lc. This difficulty arises from the failure of Bertini type theorems for foliations. In fact, to establish Theorem \ref{thm: global acc threefold} in dimension $3$, we require a gfq version of Theorem \ref{thm: global acc threefold}(2) in dimension $2$. For further details, we refer to Lemma \ref{lem: surface gloabl index with boundary}.

\medskip

\noindent\textbf{Global index theorem for foliated surfaces}. The second key ingredient in the proof of Theorem \ref{thm: global acc threefold} is the global index theorem for foliated surfaces. We prove the following result, which includes Theorem \ref{thm: global acc threefold}(2) when $\dim X=2$ but also includes the case when $B=0$:

\begin{thm}[Global index theorem for foliated surfaces]\label{thm: surface index theorem global intro}
Let $\Ii\subset [0,1]\cap\mathbb Q$ be a DCC set. Then there exists a positive integer $I$ depending only on $\Ii$ satisfying the following. Assume that $(X,\Ff,B)$ is a projective lc foliated triple such that $\dim X=2,B\in\Ii$, and $K_{\Ff}+B\equiv 0$. Then $I(K_{\Ff}+B)\sim 0$.
\end{thm}

We expect Theorem \ref{thm: surface index theorem global intro} to hold in all dimensions.

\medskip

\noindent\textbf{Functional divisors on foliations}. With Theorems \ref{thm: cbf ftriple} and \ref{thm: surface index theorem global intro} established, we are able to prove Theorem \ref{thm: global acc threefold}, which assumes $\Ii\subset\mathbb Q$. However, we still aim to obtain a complete version of the global ACC for foliated threefolds, which would hold for any DCC set $\Ii\subset [0,1]$:
\begin{conj}\label{conjthm: global acc threefold foliation}
Let $\Ii\subset [0,1]$ be a DCC set. Then there exists a finite set $\Ii_0\subset\Ii$ depending only on $\Ii$ satisfying the following. Assume that $(X,\Ff,B)$ is a projective lc foliated triple such that $\dim X\leq 3$, $K_{\Ff}+B\equiv 0$, and $B\in\Ii$. Then the coefficients of $B$ belong to $\Ii_0$.
\end{conj}
To prove Conjecture \ref{conjthm: global acc threefold foliation}, it is necessary to establish the theory of \emph{uniform rational polytopes} for \emph{functional divisors} in the context of foliations. Essentially, we aim to examine the behavior of foliated triples of the form $(X,\Ff,B(\bm{v}))$, where $\bm{v}\rightarrow B(\bm{v}): \mathbb R^m\rightarrow\Weil_{\mathbb R}(X)$ is a (linear) function. The theory of uniform rational polytopes for functional divisors suggests that the behavior of $(X,\Ff,B(\bm{v}))$ should be similar to each other when $\bm{v}$ is close to an irrational vector. More precisely, we expect the following result:

\begin{conj}\label{conjthm: uniform rational polytope foliation intro}
Let $v_1^0,\dots,v_m^0$ be positive integers and $\bm{v}_0:=(v_1^0,\dots,v_m^0)$. Then there exists an open set $U\ni \bm{v}_0$ of the rational envelope of $\bm{v}_0$ satisfying the following. 

Let $(X,\Ff,B(\bm{v}_0):=\sum_{i=1}^mv_i^0B_i)$ be any lc foliated triple of dimension $\leq 3$, where $B_i\geq 0$ are distinct Weil divisors. Then  $(X,\Ff,B(\bm{v}):=\sum_{i=1}^mv_iB_i)$ is lc for any $\bm{v}:=(v_1,\dots,v_m)\in U$.
\end{conj}

We plan to prove Conjectures \ref{conjthm: global acc threefold foliation} and \ref{conjthm: uniform rational polytope foliation intro} in upcoming work. Conjecture \ref{conjthm: uniform rational polytope foliation intro} requires different technical methods and will be treated in a separate work. However, we do provide a simplified version of Conjecture \ref{conjthm: uniform rational polytope foliation intro} in this paper:

\begin{thm}\label{thm: rational polytope foliation intro}
Let $v_1^0,\dots,v_m^0$ be positive integers, $\bm{v}_0:=(v_1^0,\dots,v_m^0)$, and $(X,\Ff,B=\sum_{i=1}^mv_i^0B_i)$ an lc foliated triple of dimension $\leq 3$, where $B_i\geq 0$ are distinct Weil divisors. Then there exists an open set $U\ni \bm{v}_0$ of the rational envelope of $\bm{v}_0$, such that $(X,\Ff,B=\sum_{i=1}^mv_iB_i)$ is lc for any $(v_1,\dots,v_m)\in U$.
\end{thm}
We note that Theorem \ref{thm: rational polytope foliation intro} is a simplified version of Conjecture \ref{conjthm: uniform rational polytope foliation intro}, where the main difference lies in the dependence of the open set $U$. In Theorem \ref{thm: rational polytope foliation intro}, $U$ may depend on $(X,\Ff,B)$, while in Conjecture \ref{conjthm: uniform rational polytope foliation intro}, $U$ only depends on $\bm{v}_0$.

As an application of Theorem \ref{thm: rational polytope foliation intro}, we can prove the following result:

\begin{thm}\label{thm: r coefficient nt abundance intro}
Let $(X,\Ff,B)$ be a projective lc foliated triple such that $\dim X\leq 3$ and $K_{\Ff}+B\equiv 0$. Then $K_{\Ff}+B\sim_{\mathbb R}0$.
\end{thm}
When $B$ is a $\mathbb{Q}$-divisor, the result of Theorem \ref{thm: r coefficient nt abundance intro} is established in \cite[Theorem 1.7]{CS21} for the case when $\dim X=3$ and $\rk\Ff=2$, in \cite[Lemma 12.2]{CS21} for the case when $\dim X=2$ and $\rk\Ff=1$, and in \cite[Theorem 1.4]{CS20} for the case when $\dim X=3$ and $\rk\Ff=1$.

We also establish a stronger version of Theorem \ref{thm: r coefficient nt abundance intro}:

\begin{thm}\label{thm: abundance num dim 0}
Let $(X,\Ff,B)$ be a projective lc foliated triple such that $\dim X\leq 3$ and $\kappa_{\sigma}(K_{\Ff}+B)=0$. Then $(X,\Ff,B)$ has a good minimal model. In particular, $\kappa_{\iota}(K_{\Ff}+B)=0$.
\end{thm}

\noindent\textit{Structure of the paper}. In Section \ref{sec: Preliminaries}, we provide a brief overview of the fundamental concepts and basic results of foliations, and introduce the definitions of gfqs and their singularities. In Section \ref{sec: algebraically integrable foliations}, we study the basic behavior of algebraically integrable foliations. In Section \ref{sec: foliated log resolution}, we introduce foliated \emph{log} resolutions and prove Theorems \ref{thm: rational polytope foliation intro} and \ref{thm: r coefficient nt abundance intro}. In Section \ref{sec: mmp foliation}, we summarize the known results of the foliated minimal model program and prove Theorem \ref{thm: abundance num dim 0}. In Section \ref{sec: Canonical bundle formula}, we establish the canonical bundle formula for foliated triples. In Section \ref{sec: Surface index theorem, global version}, we prove Theorem \ref{thm: surface index theorem global intro}. Finally, in Section \ref{sec: Proof of Theorem}, we present the proof of our main result, Theorem \ref{thm: global acc threefold}.

\medskip

\noindent\textbf{Postscript remark}. Conjectures \ref{conjthm: global acc threefold foliation} and \ref{conjthm: uniform rational polytope foliation intro} were later resolved in dimension $2$ in \cite{LMX23a} and in dimension $3$ in \cite{LMX23b}.

\medskip

\noindent\textbf{Acknowledgements}. We thank Paolo Cascini, Guodu Chen, Yen-An Chen, Omprokash Das, Christopher D. Hacon, Jingjun Han, Junpeng Jiao, Yuchen Liu, Lingyao Xie, Chenyang Xu, and Qingyuan Xue for helpful discussions.  We would like to acknowledge the assistance of ChatGPT in polishing the wording. We thank the referee for many valuable suggestions.

\section{Preliminaries}\label{sec: Preliminaries}

We will work over the field of complex numbers $\Cc$. Throughout the paper, we will mainly work with normal quasi-projective varieties to ensure consistency with the references. However, most results should also hold for normal varieties that are not necessarily quasi-projective. Similarly, most results in our paper should hold for any algebraically closed field of characteristic zero. We will adopt the standard notations and definitions in \cite{KM98, BCHM10} and use them freely. For foliations, we will generally follow the notations and definitions in \cite{CS20,ACSS21,CS21}, but there may be minor differences. For generalized pairs, we will follow the notations and definitions in \cite{HL23}.

\subsection{Sets}

\begin{defn}\label{defn: DCC and ACC}
Let $\Ii\subset\Rr$ be a set. We say that $\Ii$ satisfies the \emph{descending chain condition} (DCC) if any decreasing sequence in $\Ii$ stabilizes, and $\Ii$ satisfies the \emph{ascending chain condition} (ACC) if any increasing sequence in $\Ii$ stabilizes. 
\end{defn}

\subsection{Foliations}

\begin{defn}[Foliations, {cf. \cite[Section 2.1]{CS21}}]\label{defn: foliation}
Let $X$ be a normal variety. A \emph{foliation} on $X$ is a coherent sheaf $\Ff\subset T_X$ such that
\begin{enumerate}
    \item $\Ff$ is saturated in $T_X$, i.e. $T_X/\Ff$ is torsion free, and
    \item $\Ff$ is closed under the Lie bracket.
\end{enumerate}
The \emph{rank} of the foliation $\Ff$ is the rank of $\Ff$ as a sheaf and is denoted by $\rk\Ff$. The \emph{co-rank} of $\Ff$ is $\dim X-\rk\Ff$. The \emph{canonical divisor} of $\Ff$ is a divisor $K_\Ff$ such that $\mathcal{O}_X(-K_{\mathcal{F}})\cong\mathrm{det}(\Ff)$. We define $N_{\Ff}:=(T_X/\Ff)^{\vee\vee}$ and $N_{\Ff}^*:=N_{\Ff}^{\vee}$.

If $\Ff=0$, then we say that $\Ff$ is a \emph{foliation by points}.
\end{defn}

\begin{defn}[Singular locus]
     Let $X$ be a normal variety and $\Ff$ a rank $r$ foliation on $X$. We can associate to $\Ff$ a morphism $$\phi: \Omega_X^{[r]}\to \mathcal{O}_X(K_{\Ff})$$ defined by taking the double dual of the $r$-wedge product of the map $\Omega^1_X\to \Ff^*$, induced by the inclusion $\Ff\to T_X$. This yields a map $$\phi': (\Omega_X^{[r]}\otimes\mathcal{O}_X(-K_{\Ff}))^{\vee\vee}\to \mathcal{O}_X$$ and we define the singular locus, denoted as $\Sing~\Ff$, to be the co-support of the image of $\phi'$.
\end{defn}

\begin{defn}[Pullbacks and pushforwards, {cf. \cite[3.1]{ACSS21}}]\label{defn: pullback}
Let $X$ be a normal variety, $\Ff$ a foliation on $X$, $f: Y\dashrightarrow X$ a dominant map, and $g: X\dashrightarrow X'$ a birational map. We denote $f^{-1}\Ff$ the \emph{pullback} of $\Ff$ on $Y$ as constructed in \cite[3.2]{Dru21}. We also say that $f^{-1}\Ff$ is the \emph{induced foliation} of $\Ff$ on $Y$. We define the \emph{pushforward} of $\Ff$ on $X'$ as $(g^{-1})^{-1}\Ff$ and denote it by $g_*\Ff$.
\end{defn}

\begin{defn}[Algebraically integrable foliations, {cf. \cite[3.1]{ACSS21}}]\label{defn: algebraically integrable}
Let $X$ be a normal variety and $\Ff$ a foliation on $X$. We say that $\Ff$ is an \emph{algebraically integrable foliation} if there exists a dominant map $f: X\dashrightarrow Y$ such that $\Ff=f^{-1}\Ff_Y$, where $\Ff_Y$ is a foliation by points. In this case, we say that $\Ff$ is \emph{induced by $f$}.
\end{defn}

\begin{defn}[Invariant subvarieties, {cf. \cite[3.1]{ACSS21}}]\label{defn: f-invariant}
Let $X$ be a normal variety, $\Ff$ a foliation on $X$, and $S\subset X$ a subvariety. We say that $S$ is $\Ff$-invariant if and only if for any open subset $U\subset X$ and any section $\partial\in H^0(U,\Ff)$, we have $$\partial(\mathcal{I}_{S\cap U})\subset \mathcal{I}_{S\cap U}$$ 
where $\mathcal{I}_{S\cap U}$ is the ideal sheaf of $S\cap U$. 
\end{defn}

\begin{defn}[Non-dicritical singularities, {cf. \cite[Definition 2.10]{CS21}}]\label{defn: non-dicritical}
Let $X$ be a normal variety and $\Ff$ a foliation of rank $1$ (resp. co-rank $1$) on $X$. We say that $\Ff$ has \emph{non-dicritical} singularities if for any point (resp. closed point) $x\in X$ and any birational morphism $f: X'\rightarrow X$ such that $f^{-1}(\overline{\{x\}})$ is a divisor, each component of $f^{-1}(\overline{\{x\}})$ is $f^{-1}\Ff$-invariant.
\end{defn}

\begin{defn}[Tangent and transverse]\label{defn: tangent to foliation}
We will use two different types of tangency and transversality for foliations in this paper. Let $X$ be a normal variety, $\Ff$ a foliation on $X$, and $V\subset X$ a subvariety.
\begin{enumerate}
    \item ({cf. \cite[Definition 2.12]{CS21}}) Suppose that $\Ff$ is a co-rank $1$ foliation on $X$ with non-dicritical singularities. We say that $V$ is \emph{tangent} to $\Ff$ if for any birational morphism $f: Y\to X$ and any divisor $E$ on $Y$ such that $E$ dominates $V$, $E$ is $f^{-1}\Ff$-invariant. We say that $V$ is \emph{transverse} to $\Ff$ if $V$ is not tangent to $\Ff$.
    \item ({cf. \cite[Section 3.4]{ACSS21}}) Suppose that $\Ff$ is a rank $1$ foliation induced by a dominant rational map $X\dashrightarrow Z$. We say that $V$ is \emph{tangent} to $\Ff$ if there exists a birational morphism $\mu: X'\rightarrow X$, an equidimensional contraction $f': X'\rightarrow Z$, and a subvariety $V'\subset X'$, such that
    \begin{enumerate}
    \item $\mu^{-1}\Ff$ is induced by $f'$, and
        \item $V'$ is contained in a fiber of $f'$ and $\mu(V')=V$.
    \end{enumerate}
    We say that $V$ is \emph{transverse} to $\Ff$ if $V$ is not tangent to $\Ff$.
\end{enumerate}
(We remark that there is no conflict between these two definitions: to see this, we only need to check these two definitions coincide for rank $1$ algebraically integrable foliations on surfaces. The cases when $\dim V=0$ or $2$ are clear. When $\dim V=1$, both (1) and (2) are equivalent to the fact that $V$ is an $\Ff$-invariant curve.)

For any point $x\in V$, we say that $V$ is \emph{transverse} to $\Ff$ at $x$ if $x\not\in\Sing(X)\cup\Sing(\Ff)\cup\Sing(V)$, and for any analytic neighborhood $U$ of $x$, $T_{V}|_U\rightarrow T_X|_U$ does not factor through $T_{\Ff}|_U$. We say that $V$ is \emph{everywhere transverse} to $\Ff$ if $V$ is transverse to $\Ff$ at $x$ for any $x\in V$ (in particular, $V$ is smooth and $V$ does not intersect $\Sing(X)$ or $\Sing(\Ff)$). We say that $V$ is \emph{generically transverse} to $\Ff$ if $V$ is transverse to $\Ff$ at the generic point $\eta_V$ of $V$.
\end{defn}

\begin{defn}[Special divisors on foliations, cf. {\cite[Definition 2.2]{CS21}}]\label{defn: special divisors on foliations}
Let $X$ be a normal variety and $\Ff$ a foliation on $X$. For any prime divisor $C$ on $X$, we define $\epsilon_{\Ff}(C):=1$ if $C$ is not $\Ff$-invariant, and  $\epsilon_{\Ff}(C):=0$ if $C$ is $\Ff$-invariant. If $\Ff$ is clear from the context, then we may use $\epsilon(C)$ instead of $\epsilon_{\Ff}(C)$. For any $\Rr$-divisor $D$ on $X$, we define $$D^{\Ff}:=\sum_{C\mid C\text{ is a component of }D}\epsilon_{\Ff}(C)C.$$
Let $E$ be a prime divisor over $X$ and $f: Y\rightarrow X$ a projective birational morphism such that $E$ is on $Y$. We define $\epsilon_{\Ff}(E):=\epsilon_{f^{-1}\Ff}(E)$. It is clear that $\epsilon_{\Ff}(E)$ is independent of the choice of $f$.
\end{defn}

\subsection{Polarized foliations}

\begin{defn}
    Let $X$ be a normal variety, and $B$ and $B'$ two $\Rr$-divisors on $X$. For any prime divisor $D$ on $X$, we define $\mult_DB$ to be the coefficient of $D$ in $B$. We define
    $$B\wedge B':=\sum_{D\text{ is a prime divisor on X}}\min\{\mult_DB,\mult_DB'\}D.$$
\end{defn}

\begin{defn}[$\bb$-divisors]\label{defn: b divisors} Let $X$ be a normal quasi-projective variety. We call $Y$ a \emph{birational model} over $X$ if there exists a projective birational morphism $Y\to X$. 

Let $X\dashrightarrow X'$ be a birational map. For any valuation $\nu$ over $X$, we define $\nu_{X'}$ to be the center of $\nu$ on $X'$. A \emph{$\bb$-divisor} $\Dd$ over $X$ is a formal sum $\Dd=\sum_{\nu} r_{\nu}\nu$ where $\nu$ are valuations over $X$ and $r_{\nu}\in\mathbb R$, such that $\nu_X$ is not a divisor except for finitely many $\nu$. If in addition, $r_{\nu}\in\Qq$ for every $\nu$, then $\Dd$ is called a \emph{$\Qq$-$\bb$-divisor}. The \emph{trace} of $\Dd$ on $X'$ is the $\Rr$-divisor
$$\Dd_{X'}:=\sum_{\nu_{X'}\text{ is a divisor}}r_\nu\nu_{X'}.$$
If $\Dd_{X'}$ is $\Rr$-Cartier and $\Dd_{Y}$ is the pullback of $\Dd_{X'}$ on $Y$ for any birational model $Y$ over $X'$, we say that $\Dd$ \emph{descends} to $X'$ and $\Dd$ is the \emph{closure} of $\Dd_{X'}$, and write $\Dd=\overline{\Dd_{X'}}$. 

Let $X\rightarrow U$ be a projective morphism and assume that $\Dd$ is a $\bb$-divisor over $X$ such that $\Dd$ descends to some birational model $Y$ over $X$. If $\Dd_Y$ is nef$/U$ (resp. base-point-free$/U$, semi-ample$/U$), then we say that $\Dd$ is \emph{nef}$/U$ (resp. \emph{base-point-free}$/U$, \emph{semi-ample}$/U$). If $\Dd_Y$ is a Cartier divisor, then we say that $\Dd$ is \emph{$\bb$-Cartier}. If $\Dd_Y$ is a $\Qq$-Cartier $\Qq$-divisor, then we say that $\Dd$ is \emph{$\Qq$-$\bb$-Cartier}. If $\Dd$ can be written as an $\Rr_{\geq 0}$-linear combination of nef$/U$ $\bb$-Cartier $\bb$-divisors, then we say that $\Dd$ is \emph{NQC}$/U$.

Let $X\rightarrow U$ be a projective morphism and assume that $\Dd$ and $\Dd'$ are two $\bb$-divisors over $X$. We write $\Dd\sim_{\mathbb R,U}\Dd'$ (resp. $\Dd\sim_{\mathbb Q,U}\Dd'$) if there exists a birational model $Y$ of $X$ such that $\Dd$ and $\Dd'$ descend to $Y$, and $\Dd_Y\sim_{\mathbb R,U}\Dd'_Y$ (resp. $\Dd_Y\sim_{\mathbb Q,U}\Dd'_Y$). 

We let $\bm{0}$ be the $\bb$-divisor $\bar{0}$.
\end{defn}

\begin{defn}[Generalized foliated quadruples]
A \emph{generalized foliated sub-quadruple} (\emph{sub-gfq} for short) $(X,\Ff,B,\Mm)/U$ consists of a normal quasi-projective variety $X$, a foliation $\Ff$ on $X$, an $\Rr$-divisor $B$ on $X$, a projective morphism $X\rightarrow U$, and an NQC$/U$ $\bb$-divisor $\Mm$ over $X$, such that $K_{\Ff}+B+\Mm_X$ is $\mathbb R$-Cartier. If $B\geq 0$, then we say that $(X,\Ff,B,\Mm)/U$ is a \emph{generalized foliated quadruple} (\emph{gfq} for short). If $U=\{pt\}$, we usually drop $U$ and say that $(X,\Ff,B,\Mm)$ is \emph{projective}. 

Let $(X,\Ff,B,\Mm)/U$ be a (sub-)gfq. If $\Mm=\bm{0}$, then we may denote $(X,\Ff,B,\Mm)/U$ by $(X,\Ff,B)/U$ or $(\Ff,B)$, and say that $(X,\Ff,B)$ is a \emph{foliated (sub-)triple} (\emph{f-(sub-)triple} for short) and $(\Ff,B)$ is a \emph{foliated (sub-)pair} (\emph{f-(sub-)pair} for short). If $\Ff=T_X$, then we may denote $(X,\Ff,B,\Mm)/U$ by $(X,B,\Mm)/U$, and say that $(X,B,\Mm)/U$ is a \emph{generalized (sub-)pair} (\emph{g-(sub-)pair} for short). If $\Mm=\bm{0}$ and $\Ff=T_X$, then we may denote $(X,\Ff,B,\Mm)/U$ by $(X,B)/U$, and say that $(X,B)/U$ is a \emph{(sub-)pair}. 

A (sub-)gfq (resp. f-(sub-)triple, f-(sub-)pair, g-(sub-)pair, (sub-)pair) $(X,\Ff,B,\Mm)/U$ (resp. $(X,\Ff,B)/U$, $(\Ff,B)/U$, $(X,B,\Mm)/U$, $(X,B)/U$) is called a \emph{$\mathbb Q$-(sub-)gfq} (resp. \emph{$\mathbb Q$-f-(sub-)triple, $\mathbb Q$-f-(sub-)pair, $\mathbb Q$-g-(sub-)pair, $\mathbb Q$-(sub-)pair} if $B$ is a $\mathbb Q$-divisor and $\Mm$ is a $\mathbb Q$-$\bb$-divisor.

A \emph{foliated germ} (\emph{f-germ} for short) $(X\ni x,\Ff,B)$ consists of a foliated triple $(X,\Ff,B)/U$ and a (not necessarily closed) point $x\in X$ where the projective morphism $X\rightarrow U$ is the identity morphism (in particular, $U=X$). We also say that $(\Ff\ni x,B)$ is an \emph{f-germ}. If $\Ff=T_X$, then we say that $(X\ni x,B)$ is a \emph{germ}. 
\end{defn}

\begin{nota}
In the previous definition, if $U$ is not important, we may also drop $U$. This usually happens when we emphasize the structures of $(X,\Ff,B,\Mm)$ which are independent of the choice of $U$, such as the singularities of $(X,\Ff,B,\Mm)$. In addition, if $B=0$, then we may drop $B$.
\end{nota}

\begin{defn}[Singularities of gfqs]\label{defn: gfq singularity}
Let $(X,\Ff,B,\Mm)$ be a (sub-)gfq. For any prime divisor $E$ over $X$, let $f: Y\rightarrow X$ be a birational morphism such that $E$ is on $Y$, and suppose that
$$K_{\Ff_Y}+B_Y+\Mm_Y:=f^*(K_\Ff+B+\Mm_X)$$
where $\Ff_Y:=f^{-1}\Ff$. We define $a(E,\Ff,B,\Mm):=-\mult_EB_Y$ to be the \emph{discrepancy} of $E$ with respect to $(X,\Ff,B,\Mm)$. It is clear that $a(E,\Ff,B,\Mm)$ is independent of the choice of $Y$. If $\Mm=\bm{0}$, then we let $a(E,\Ff,B):=a(E,\Ff,B,\Mm)$. If $\Ff=T_X$, then we let $a(E,X,B,\Mm):=a(E,\Ff,B,\Mm)$. If $\Mm=\bm{0}$ and $\Ff=T_X$, then we let $a(E,X,B):=a(E,\Ff,B,\Mm)$.

We say that $(X,\Ff,B,\Mm)$ is \emph{(sub-)lc} (resp. \emph{(sub-)klt}) if $a(E,\Ff,B,\Mm)\geq -\epsilon_{\Ff}(E)$ (resp. $>-\epsilon_{\Ff}(E)$) for any prime divisor $E$ over $X$. We say that $(X,\Ff,B,\Mm)$ is \emph{(sub-)canonical} (resp. \emph{(sub-)terminal}) if $a(E,\Ff,B,\Mm)\geq 0$ (resp. $>0$) for any prime divisor $E$ that is exceptional over $X$.

Let $(X\ni x,\Ff,B)$ be an f-germ. We say that $(X\ni x,\Ff,B)$ and $(\Ff\ni x,B)$ are \emph{lc} (resp. \emph{klt, canonical, terminal}) if $a(E,\Ff,B)\geq-\epsilon(E)$ (resp. $>-\epsilon(E),\geq 0,>0$) for any prime divisor $E$ over $X\ni x$.
\end{defn}

\begin{rem}
In the literature of generalized pairs, ``generalized lc (glc), generalized klt (gklt)" are often used, while in the literature of foliations, ``foliated lc, foliated klt" are often used. Since we introduce the concept of ``gfq", we feel that it is a good practice to simply use ``lc, klt", removing the adjective ``generalized" or ``foliated", even if we only care about the case when $\Ff=T_X$ (the usual generalized pairs) or the case when $\Mm=\bm{0}$ (the usual foliated triples). Similarly, we shall use ``dlt" instead of ``gdlt" or ``F-dlt". 
\end{rem}

\subsection{Simple singularities}

\begin{defn}[Non-resonant]
Let $k$ be a positive integer and $\lambda_1,\dots,\lambda_k$ real numbers. We say that $\lambda_1,\dots,\lambda_k$ are \emph{non-resonant} if for any non-negative integers $a_1,\dots,a_k$ such that $\sum_{i=1}^ka_i\lambda_i=0$, we have $a_i=0$ for all $i$.
\end{defn}

\begin{defn}
We will use two types of \emph{simple singularities} in this paper. Let $d$ be a positive integer, $X$ a normal variety of dimension $d$, and $\Ff$ a foliation on $X$.
\begin{enumerate}
    \item (\cite[Definition 2.8]{CS21}) Suppose that $\rk\Ff=\dim X-1$. For any point $o\in X$, we say that $o$ is a \emph{simple singularity} of $\Ff$ if $X$ is smooth near $o$, and there exist formal local coordinates $x_1,\dots,x_d$ near $o$ and some integer $1\leq r\leq d$, such that $N_{\Ff}^*$ is generated by one of the following two types of $1$-forms $\omega$:
    \begin{enumerate}
        \item $$\omega=x_1\dots x_r\sum_{i=1}^r\lambda_i\frac{d x_i}{x_i}$$
        where $\lambda_1,\dots,\lambda_r$ are non-resonant.
        \item
        There is an integer $k \le r$ such that
        $$w=x_1\dots x_r\left(\sum_{i=1}^kp_i\frac{dx_i}{x_i}+\psi(x_1^{p_1}\dots x_k^{p_k})\sum_{i=2}^r\lambda_i\frac{dx_i}{x_i}\right)$$
    \end{enumerate}
    where $p_1,\dots,p_k$ are positive integers such that $\gcd(p_1,\dots,p_k)=1$, $\psi$ is a formal power series that is not a unit, and $\lambda_2,\dots,\lambda_k$ are non-resonant.
    \item (\cite[Definition 2.32]{CS20}) Suppose that $d=3$ and $\rk\Ff=1$. For any point $o\in X$, we say that $o$ is a \emph{simple singularity} of $\Ff$ if one of the following holds.
    \begin{enumerate}
    \item $\Ff$ is terminal near $o$ and no component of $\Sing(X)$ passing through $o$ is $\Ff$-invariant.
    \item   $(X\ni o,\Ff)$ is formally isomoprhic to $(X_0\ni x_0,\Ff_0)$, where $(X_0\ni x_0):=(\mathbb C^3\ni 0)/\bm{\mu}$ where $\bm{\mu}$ is the $\mathbb Z_2$-action given by $(x,y,z)\rightarrow (y,x,-z)$, and $\Ff_0$ is the rank $1$ foliation on $X_0$ induced by the vector field
    $$\partial:=\left(\left(1+a(xy,z)\right)x\frac{\partial}{\partial x}-\left(1+a(xy,-z)\right)y\frac{\partial}{\partial y}+c(xy,z)\frac{\partial}{\partial z}\right)/\bm{\mu}$$
    where $a$ and $c$ are formal functions in two variables such that $c$ is not a unit and $c(xy,z)=c(xy,-z)$ for any $x,y,z$. 
    \item $\Ff$ has canonical singularities and $X$ is smooth at $o$.
    \end{enumerate}
\end{enumerate}
We say that $\Ff$ has \emph{at most simple singularities} if $o$ is a simple singularity of $\Ff$ for any $o\in\Sing(\Ff)$.
\end{defn}

\begin{lem}\label{lem: simple singularity nondicritical}
Let $X$ be a normal variety and $\Ff$ a foliation on $X$. Suppose that either $\rk\Ff=\dim X-1$, or $\dim X=3$ and $\rk\Ff=1$. If $\Ff$ has canonical singularities, then $\Ff$ is non-dicritical.
\end{lem}
\begin{proof}
If $\rk\Ff=\dim X-1$, then the lemma follows from \cite[Remark 2.13]{CS21}. If $\dim X=3$ and $\rk\Ff=1$, then the lemma follows from \cite[Corollary III.i.4]{MP13}.
\end{proof}

\begin{lem}[{cf. \cite[Lemma 2.34]{CS20}}]\label{lem: simple singularities cyc}
Let $X$ be a normal variety and $\Ff$ a foliation on $X$ such that $\dim X=3$ and $\rk\Ff=1$. If $\Ff$ has at most simple singularities, then $X$ has at most cyclic quotient singularities. In particular, $X$ is $\Qq$-factorial klt.
\end{lem}

\section{Algebraically integrable foliations}\label{sec: algebraically integrable foliations}

The goal of this section is to prove the following theorem:

\begin{thm}\label{thm: subfoliation algebraic integrable}
Let $\Ff$ be a foliation on a normal projective variety $X$ such that $K_{\Ff}$ is not pseudo-effective. Then there exists an algebraically integrable foliation $\Ee$ such that $0\not=\Ee\subset\Ff$.
\end{thm}

Theorem \ref{thm: subfoliation algebraic integrable} is essentially a corollary of \cite[Theorem 1.1]{CP19} and should be well-known to experts. Since we cannot find any direct reference, we give a proof in this paper.

\begin{defn}[Slope, cf. {\cite[Definitions 2.9 and 2.20]{GKP16},~\cite[Section 2.1]{Ou23}}]\label{defn: slope}
Let $X$ be a $\mathbb Q$-factorial normal projective variety, $\alpha$ a movable class on $X$, and $\Ff\not=0$ a torsion-free coherent sheaf on $X$. We define
$$\mu_{\alpha}(\Ff):=\frac{c_1(\Ff)\cdot\alpha}{\rk\Ff}$$
to be the \emph{slope} of $\Ff$ with respect to $\alpha$. We define
$$\mu_{\alpha}^{\max}(\Ff):=\sup\{\mu_{\alpha}(\Gg)\mid 0\not=\Gg\subset\Ff\text{ is a coherent subsheaf}\},$$
and
$$\mu_{\alpha}^{\min}(\Ff):=\inf\{\mu_{\alpha}(\mathcal{Q})\mid \mathcal{Q}\text{ is a non-zero torsion-free quotient sheaf of } \mathcal{\Ff}\}.$$

These definitions coincide with \cite[Definition 2.3]{CP19} when $X$ is smooth.
\end{defn}

\begin{defn}[$\alpha$-semistable, {\cite[Definition 2.10]{GKP16}}] Let $X$ be a $\mathbb Q$-factorial normal projective variety, $\alpha$ a movable class on $X$, and $\Ff\not=0$ a torsion-free coherent sheaf on $X$. We say that $\Ff$ is \emph{$\alpha$-semistable} if $\mu_{\alpha}(\Gg)\leq\mu_{\alpha}(\Ff)$ for any coherent subsheaf $0\not=\Gg\subset\Ff$. A \emph{maximal $\alpha$-destabilizing sheaf} of $\Ff$ is an $\alpha$-semistable sheaf $0\not=\Gg$ saturated in $\Ff$ such that
\begin{enumerate}
    \item $\mu_{\alpha}^{\max}(\Ff)=\mu_{\alpha}(\Gg)$, and
    \item $\Gg'\subset\Gg$ for any subsheaf $\Gg'\subset\Ff$ such that $\mu_{\alpha}^{\max}(\Ff)=\mu_{\alpha}(\Gg')$.
\end{enumerate}
\end{defn}

\begin{lem}[{cf. \cite[Lemma~4.10]{CP19}}]\label{lem: a foliation criterion}
Let $X$ be a smooth projective variety, $\alpha$ a movable class on $X$, and $\Gg\subset\Ff\subset T_X$ coherent sheaves. Assume that $\Ff$ is a foliation, $\Gg$ is saturated in $\Ff$, $\mu_{\alpha}^{\min}(\Gg)>0$, and
$$2\mu_{\alpha}^{\min}(\Gg)>\mu_{\alpha}^{\max}(\Ff/\Gg).$$
Then $\Gg$ is a foliation.
\end{lem}
\begin{proof}
Since $\Gg$ is saturated in $\Ff$ and $\Ff$ is saturated in $T_X$, $\Ff/\Gg$ and $T_X/\Ff$ are torsion-free. By the exact sequence
$$0\rightarrow\Ff/\Gg\rightarrow T_X/\Gg\rightarrow T_X/\Ff\rightarrow0,$$
we deduce that $T_X/\Gg$ is torsion-free and thus $\Gg$ is saturated in $T_X$. The proof for the fact that $\Gg$ is closed under the Lie bracket is identical to \cite[Lemma 4.10]{CP19}, so we omit it (although \cite[Lemma 4.10]{CP19} requires that $\Gg$ and $\Ff$ are distributions).
\end{proof}

\begin{lem}\label{lem: algebraically integral when chern class none pe}
Let $\Ff$ be a foliation on a smooth projective variety $X$ such that $K_{\Ff}$ is not pseudo-effective. Then there exists an algebraically integrable foliation $\Ee$ such that $0\not=\Ee\subset\Ff$.
\end{lem}

\begin{proof}
Since $K_{\Ff}$ is not pseudo-effective, by \cite[Theorem 0.2]{BDPP13}, there exists a movable class $\alpha$ on $X$ such that $\alpha\cdot K_\Ff<0$. Thus $$\mu_\alpha(\Ff)=\frac{\alpha\cdot c_1(\Ff)}{\rank\Ff}=\frac{\alpha\cdot(-K_\Ff)}{\rank\Ff}>0.$$
By \cite[Corollary 2.26]{GKP16}, there exist a positive integer $r$ and a unique Harder-Narasimhan filtration $$0=\Ee_0\subsetneq \Ee_1\subsetneq\cdots\subsetneq \Ee_r=\Ff,$$
i.e. a filtration such that $\Ee_k/\Ee_{k-1}$ is torsion-free and $\alpha$-semistable for any $1\leq k\leq r$, and $\mu_{\alpha}(\Ee_k/\Ee_{k-1})<\mu_{\alpha}(\Ee_{k-1}/\Ee_{k-2})$ for any $2\leq k\leq r$. For any $1\leq k\leq r$, $\Ee_k/\Ee_{k-1}$ is the maximal $\alpha$-destabilizing sheaf of $\Ff/\Ee_{k-1}$, hence $\mu_{\alpha}(\Ee_k/\Ee_{k-1})=\mu_{\alpha}^{\max}(\Ff/\Ee_{k-1})$. We have
 $$\mu_{\alpha}^\mathrm{min}(\Ee_1)=\mu_{\alpha}(\Ee_1)=\mu_{\alpha}^\mathrm{max}(\Ff)\geq \mu_{\alpha}(\Ff)>0.$$
 If $r=1$, then $\Ff=\Ee_1$ and the lemma follows from  \cite[Theorem 1.1]{CP19}. Thus we may assume that $r\geq 2$. Since $\mu_{\alpha}^\mathrm{min}(\Ee_1)>0$, $$2\mu_{\alpha}^\mathrm{min}(\Ee_1)>\mu_{\alpha}^\mathrm{min}(\Ee_1)=\mu_{\alpha}(\Ee_{1})>\mu_{\alpha}(
\Ee_2/\Ee_{1})=\mu_{\alpha}^\mathrm{max}(\Ff/\Ee_1).$$
Since $\Ee_1$ is the maximal $\alpha$-destabilizing sheaf of $\Ff$, $\Ee_1$ is saturated in $\Ff$. By Lemma \ref{lem: a foliation criterion}, $\Ee:=\Ee_1$ is a foliation. By \cite[Theorem 1.1]{CP19}, $\Ee$ is algebraically integrable.
\end{proof}

\begin{proof}[Proof of Theorem \ref{thm: subfoliation algebraic integrable}]
We let $f: Y\rightarrow X$ be a resolution of $X$. Since $K_{\Ff}=f_*K_{f^{-1}\Ff}$, $K_{\Ff}$ is not pseudo-effective, and the pushforward of any pseudo-effective divisor under birational morphisms is still pseudo-effective, $K_{f^{-1}\Ff}$ is not pseudo-effective. 

By Lemma \ref{lem: algebraically integral when chern class none pe}, there exists an algebraically integrable foliation $\Ee_Y$ such that $0\not=\Ee_Y\subset f^{-1}\Ff$. We let $\Ee':=f_*\Ee_Y\not=0$ (as a foliation). The foliation $\Ee'$ is also algebraically integrable. Then we have
$$\Ee'\cong(f_*\Ee_Y)^{\vee\vee}\subset(f_*(f^{-1}\Ff))^{\vee\vee}\cong\Ff$$ 
as sheaves. Thus there exists $0\not=\Ee\subset\Ff$ which is isomorphic to $\Ee'$, and $\Ee$ is algebraically integrable.
\end{proof}

\section{Foliated log resolutions and proofs of Theorems \ref{thm: rational polytope foliation intro} and \ref{thm: r coefficient nt abundance intro}}\label{sec: foliated log resolution}

It is well-known that any foliation of dimension $\leq 3$ has a \emph{foliated resolution}. Nevertheless, in practice, we need to use the log version of a foliated resolution, i.e. a \emph{foliated log resolution}. The existence of such a log resolution in dimension $\leq 3$ is actually known in literature yet not precisely summarized, so we summarize them in this section. As an application, we prove Theorem \ref{thm: rational polytope foliation intro}.

\subsection{Foliated log resolutions}

\begin{defn}[{\cite[Definition 3.1]{CS21}}]\label{defn: folaited log smooth corank 1}
Let $(X,\Ff,B)$ be an f-sub-triple such that $\rk\Ff=\dim X-1$. We say that $(X,\Ff,B)$ is \emph{foliated log smooth} if
\begin{enumerate}
    \item $(X, B)$ is log smooth,
    \item $\Ff$ has at most simple singularities, and
    \item for any closed point $x\in B^{\Ff}$, suppose that $D_1,\dots,D_k$ are the $\Ff$-invariant divisors passing through $x$, then $B^{\Ff}\cup(\cup_{i=1}^kD_i$) is a normal crossing divisor at $x$.
\end{enumerate}
\end{defn}

We need the following theorem on \emph{log} resolutions of singularities for foliations of rank $1$ in dimension $3$. First we define foliated log smoothness for foliations of rank $1$ in dimension $3$.

\begin{defn}\label{defn: foliated log smooth rank 1}
Let $(X,\Ff,B)$ be a f-sub-triple such that $\dim X=3$ and $\rk\Ff=1$. We say that $(X,\Ff,B)$ is \emph{foliated log smooth} if
\begin{enumerate}
    \item $\Supp B$ is simple normal crossing,
    \item $X$ has at most $\mathbb Z/2\mathbb Z$-quotient singularities and $\Ff$ has at most simple singularities, and
    \item $B^{\Ff}$ is everywhere transverse to $\Ff$. In particular, $B^{\Ff}$ does not intersect $\Sing(X)\cup\Sing(\Ff)$, and the components of $B^{\Ff}$ are disjoint with each other.
\end{enumerate}
Note that although we call $(X,\Ff,B)$ foliated log smooth, $X$ may not be smooth.
\end{defn}

\begin{lem}\label{lem: corank log smooth imply flc}
Let $(X,\Ff,B)$ be an f-sub-triple such that $(X,\Ff,B)$ is foliated log smooth, and either $\rk\Ff=\dim X-1$, or $\dim X=3$ and $\rk\Ff=1$. Then $(X,\Ff,B^{\Ff})$ is lc.
\end{lem}
\begin{proof}
When  $\rk\Ff=\dim X-1$, the lemma follows from \cite[Lemma 3.3]{CS21}. Thus we may assume that $\dim X=3$ and $\rk\Ff=1$. We only need to show that $(X\ni x,\Ff,B^{\Ff})$ is lc for any point $x\in X$. There is nothing left to prove when $\dim\overline{\{x\}}\geq 2$, so we may assume that $\dim\overline{\{x\}}\leq 1$. Since $\Ff$ has simple singularities, we may assume that $x\in\Supp B^{\Ff}$. Since $B^{\Ff}$ does not intersect $\Sing(X)\cup\Sing(\Ff)$, $x\not\in\Sing(X)\cup\Sing(\Ff)$.

Let $f: Y\rightarrow X$ be the ordinary blow-up of $x$ with exceptional divisor $E$. By Lemma \ref{lem: simple singularity nondicritical}, $\Ff$ is non-dicritical, so $E$ is $\Ff_Y$-invariant. Then $a(E,\Ff,0)=2-\dim\overline{\{x\}}$. Since the components of $B^{\Ff}$ are disjoint with each other, $x$ is only contained in one component of $B$, hence $\mult_EB^{\Ff}=1$. Thus $a(E,\Ff,B^{\Ff})=1-\dim\overline{\{x\}}\geq 0=\epsilon(E)$. Let $K_{\Ff_Y}+B_Y:=f^*(K_{\Ff}+B^\Ff)$, where $\Ff_Y:=f^{-1}\Ff$. Then $(Y,\Ff_Y,B_Y)$ is foliated log smooth. To see this, we note that Definition \ref{defn: foliated log smooth rank 1}(1)(2) follow immediately from our construction, and Definition \ref{defn: foliated log smooth rank 1}(3) follows from the fact that $B^{\Ff}$ is transverse to $\Ff$ at $x$ and $\Ff$ is non-dicritical. By replacing $(X,\Ff,B)$ with $(Y,\Ff_Y, B_Y)$ and repeating our arguments, we have that $a(D,\Ff,B)\geq \epsilon(D)$ for any prime divisor $D$ over $X\ni x$. The lemma follows.
\end{proof}

\begin{defn}\label{defn: log resolution}
Let $X$ be a normal variety, $\Ff$ a foliation on $X$, and $B$ an $\Rr$-divisor on $X$, such that either $\rk\Ff=\dim X-1$ or $\dim X=3$ and $\rk\Ff=1$. A \emph{foliated log resolution} of $(X,\Ff,B)$ is a projective birational morphism $f: Y\rightarrow X$ such that $(Y,\Ff_Y:=f^{-1}\Ff,B_Y:=f^{-1}_*B+\Exc(f))$ is foliated log smooth, where $\Exc(f)$ is the reduced exceptional divisor of $f$. A \emph{foliated resolution} of $\Ff$ is a foliated log resolution of $(X,\Ff,0)$.
\end{defn}

\begin{thm}\label{thm: existence of resolution}
Let $(X,\Ff,B)$ be a foliated sub-triple such that $\dim X\leq 3$ and $1\leq\rk\Ff\leq\dim X-1$. Then there exists a foliated log resolution of $(X,\Ff,B)$.
\end{thm}
\begin{proof}
 When $\dim X=2$ and $\rk\Ff=1$, the theorem essentially follows from \cite{Sei68} (we refer to \cite[Pages 908--912]{Can04} for a detailed explanation). When $\dim X=3$ and $\rk\Ff=2$, the theorem follows from \cite[Theorem]{Can04}. Thus we may assume that $\dim X=3$ and $\rk\Ff=1$.

By \cite[Scholion]{MP13}, there exists a birational projective morphism $g: W\rightarrow X$ such that 
\begin{enumerate}
\item $B_W:=g^{-1}_*B+\Exc(g)$ is simple normal crossing, where $\Exc(g)$ is the reduced exceptional divisor of $g$,
    \item $W$ has at most $\mathbb Z/2\mathbb Z$-quotient singularities and $\Ff_W:=g^{-1}\Ff$ has at most simple singularities, and
    \item for any component $D$ of $B_W$, $D$ is
    \begin{itemize}
        \item either $\Ff_W$-invariant, or
        \item everywhere transverse to $\Ff_W$.
    \end{itemize}
    (Note that the notion ``everywhere transverse" is undefined in \cite{MP13}. Nevertheless, the ``everywhere transverse" condition of \cite[Scholion]{MP13} is essentially given by the construction in \cite[III.ii.3 Log resolution]{MP13}, which coincides with Definition \ref{defn: tangent to foliation}.)
\end{enumerate}
Since $\Ff_W$ has simple singularities, by Lemma \ref{lem: simple singularity nondicritical}, $\Ff_W$ is non-dicritical. Thus we may let $h: Y\rightarrow W$ be a birational morphism which blows up the strata of $\Supp B_W$, and let $f:=g\circ h$, $\Ff_Y:=f^{-1}\Ff$, and $K_{\Ff_Y}+B_Y:=f^*(K_X+B)$. Then the components of $B_Y$ that are not $\Ff_Y$-invariant are disjoint with each other, hence $(Y,\Ff_Y,B_Y)$ is foliated log smooth.
\end{proof}

\subsection{Proofs of Theorems \ref{thm: rational polytope foliation intro} and \ref{thm: r coefficient nt abundance intro}}

\begin{lem}\label{lem: weak perturbation}
Let $c, m$ be two non-negative integers, $r_1,\dots,r_c$ real numbers such that  $1,r_1,\dots,r_c$ are $\mathbb Q$-linearly independent, and $s_1,\dots,s_m:\mathbb R^{c+1}\rightarrow\mathbb R$ $\Qq$-linear functions. Let $(X,\Ff,B)$ be an lc (resp. klt, dlt) f-triple such that $\dim X\leq 3$, $1\leq\rk\Ff\leq\dim X-1$, and $B=\sum_{i=1}^ms_i(1,\bm{r})B_i$, where $\bm{r}:=(r_1,\dots,r_c)$ and $B_i\geq 0$ are Weil divisors on $X$. We let $B(\bm{v}):=\sum_{i=1}^ms_i(1,\bm{v})B_i$ for any $\bm{v}\in\mathbb R^c$. Then there exists an open set $U\ni\bm{r}$ of $\mathbb R^c$, such that $(X,\Ff,B(\bm{v}))$ is lc for any $\bm{v}\in U$ (resp. klt, dlt).
\end{lem}
\begin{proof}
By \cite[Lemma 5.3]{HLS19}, for any $\bm{v}\in\mathbb R^c$, $K_{\Ff}+B(\bm{v})$ is $\Rr$-Cartier. By Theorem \ref{thm: existence of resolution}, there exists a foliated log resolution $f: Y\rightarrow X$ of $(X,\Ff,B)$ with $\Ff_Y:=f^{-1}\Ff$. Moreover, if $(X,\Ff,B)$ is dlt, $\dim X=3$, and $\rk\Ff=2$, then we may choose $f$ such that $a(E,\Ff,B(\bm{r}))>-\epsilon_{\Ff_Y}(E)$ for any $f$-exceptional prime divisor $E$ on $Y$. 

We let $K_{\Ff_Y}+B_Y(\bm{v}):=f^*(K_{\Ff}+B(\bm{v}))$ for any $\bm{v}\in\mathbb R^c$. Then for any prime divisor $D$ on $Y$ that is contained in $\Supp f^{-1}_*B\cup\Exc(f)$, $\mult_DB_Y(\bm{r})\leq\epsilon_{\Ff_Y}(D)$, and $\mult_DB_Y(\bm{v})$ is independent of $\bm{v}$ if $\mult_DB_Y(\bm{r})\in\mathbb Q$.

Since there are only finitely many prime divisors contained in $\Supp f^{-1}_*B\cup\Exc(f)$ and $\bm{v}\rightarrow\mult_DB_Y(\bm{v})$ is a continuous function, there exists an open set $U\ni\bm{r}$ of $\mathbb R^c$, such that for any prime divisor $D$ on $Y$ that is contained in $\Supp f^{-1}_*B\cup\Exc(f)$, $\mult_DB_Y(\bm{v})\leq\epsilon_{\Ff_Y}(D)$, and $\mult_DB_Y(\bm{v})<\epsilon_{\Ff_Y}(D)$ if $\mult_DB_Y(\bm{r})<\epsilon_{\Ff_Y}(D)$. Thus $(Y,\Ff_Y,B_Y(\bm{v}))$ is lc (resp. klt, dlt) for any $\bm{v}\in\mathbb R^c$. 
\end{proof}

\begin{proof}[Proof of Theorem \ref{thm: rational polytope foliation intro}]
When $1\leq\rk\Ff\leq\dim X-1$, the theorem follows from Lemma \ref{lem: weak perturbation}. When $\rk\Ff=0$, lc implies that $B=0$ so there is nothing left to prove. When $\rk\Ff=\dim X$, the theorem is well-known, cf. \cite[Theorem 5.6]{HLS19}.
\end{proof} 

\begin{proof}[Proof of Theorem \ref{thm: r coefficient nt abundance intro}]
When $\rk\Ff=0$, we have $K_\Ff=B=0$ and there is nothing left to prove. When $\rk\Ff=\dim X$, the theorem is well-known, cf. \cite[Theorem 6.2]{Gon11}. So we may assume that $1\leq\rk\Ff\leq\dim X-1$.

By Lemma \ref{lem: weak perturbation} and \cite[Lemma 5.3]{HLS19}, there exist real numbers $a_1,\dots,a_k\in (0,1]$ and $\Qq$-divisors $B_1,\dots,B_k$ on $X$, such that $\sum_{i=1}^ka_i=1$, $B=\sum_{i=1}^ka_iB_i$,  $(X,\Ff,B_i)$ is lc for each $i$, and $K_{\Ff}+B_i\equiv 0$ for each $i$. 

By \cite[Theorem 1.4]{CS20} (when $\rk\Ff=1$ and $\dim X=3$), or \cite[Theorem 8.1 and Lemma 12.2]{CS21} (when $\rk\Ff=1$ and $\dim X=2$), or \cite[Theorem 1.7]{CS21}  (when $\rk\Ff=2$ and $\dim X=3$), $\kappa(K_{\Ff}+B_i)=0$ for each $i$, and thus $K_{\Ff}+B_i\sim_{\mathbb Q}0$ for each $i$. Thus $K_{\Ff}+B\sim_{\mathbb R}0$.
\end{proof}

\section{Minimal model program for foliations and proof of Theorem \ref{thm: abundance num dim 0}}\label{sec: mmp foliation}

The foundations of the minimal model program for foliated surfaces and threefolds have been established in \cite{McQ08,Bru15} and \cite{CS20,Spi20,CS21,SS22} respectively. However, there are many results on rank one foliated pairs in dimension $3$ that are essentially proved in \cite{CS20}, but not yet precisely stated. This is because most parts of \cite{CS20} focus on $K_{\Ff}$-negative extremal rays but not $(K_{\Ff}+B)$-negative extremal rays, which makes those results difficult to be directly referred to. In this section, we resolve this problem. 

\begin{rem}
In this section, there are many results where the additional condition ``being projective" on a foliated variety $(X,\Ff)$ is assumed, e.g. Lemma \ref{lem: rel version 8.1 CS20}, Theorem \ref{thm: fdlt modification} and Theorem \ref{thm: can run mmp}. This condition is expected to be unnecessary as the relative versions of the theorems are still expected to hold (e.g. \cite{SS22} for the co-rank one case). Nevertheless, since the key references \cite{CS20,CS21} only deal with the projective case, we will be consistent with these references. 

We remark that, although most results for foliated MMP in dimension $\leq 3$ in the relative setting are expected to be deduced from the same lines of the proofs in the projective setting, there may be some tricky details to deal with. For example, the bend and break technique \cite[Corollary 2.28]{Spi20}, which is essentially used in the cone theorem of rank one foliations in dimension $3$ \cite[Proof of Theorem 2.36]{CS20}, no longer works straightforwardly. On the other hand, when the base variety is projective, we can reduce the question to the projective case (see Lemma \ref{lem: rel version 8.1 CS20} below). We also expect the methods in \cite{KMM87} to be applied to the relative foliated MMP even when the base variety is no longer projective.
\end{rem}

\subsection{Dlt models} We define dlt singularities and dlt models for foliations. They are defined for the co-rank one case, but not for the rank one case in dimension $3$ in literature.

\begin{defn}[Dlt singularities]\label{def: fdlt}
Let $(X,\Ff,B)$ be an f-triple. 
\begin{enumerate}
    \item Suppose that $\rk\Ff=\dim X-1$. We say that $(X,\Ff,B)$ and $(\Ff,B)$ are \emph{dlt} if 
\begin{enumerate}
\item $(X,\Ff,B)$ is lc,
\item every component of $B$ is generically transverse to $\Ff$, and
\item there exists a foliated log resolution of $(X,\Ff,B)$ which only extracts divisors $E$ of discrepancy $>-\epsilon(E)$.
\end{enumerate}
\item Suppose that $\dim X=3$ and $\rk\Ff=1$.  We say that $(X,\Ff,B)$ and $(\Ff,B)$ are \emph{dlt} if 
\begin{enumerate}
\item $(X,\Ff,B)$ is lc, and
\item $\Ff$ has simple singularities.
\end{enumerate}
\end{enumerate}
We remark that dlt implies non-dicritical (cf. \cite[Theorem 11.3]{CS21} and Lemma \ref{lem: simple singularity nondicritical}).
\end{defn}

\begin{defn}[Dlt modification]\label{defn: fdlt modification}
Let $(X,\Ff,B)$ be an lc f-triple. An \emph{dlt modification} of $(X,\Ff,B)$ is a birational morphism $f: Y\rightarrow X$ satisfying the following. Let $\Ff_Y:=f^{-1}\Ff$, $E$ the reduced exceptional divisor of $f$, and  $B_Y:=f^{-1}_*B+E^{\Ff_Y}$.
\begin{enumerate}
    \item $Y$ is $\mathbb Q$-factorial klt, 
    \item $K_{\Ff_Y}+B_Y=f^*(K_{\Ff}+B)$, and
    \item $(Y,\Ff_Y,B_Y)$ is dlt.
\end{enumerate}
We say that $(Y,\Ff_Y,B_Y)$ is a \emph{dlt model} of $(X,\Ff,B)$ and $(\Ff_Y,B_Y)$ is a \emph{dlt model} of $(\Ff,B)$.
\end{defn}

\begin{rem}
We define ``dlt$=$lc$+$simple singularities" for rank one foliations in dimension $3$ because it satisfies the following three natural properties of usual dlt pairs (with property (1) for $\Qq$-factorial dlt pairs): 
\begin{enumerate}
    \item the ambient variety is klt,
    \item log smoothness implies dlt, and
    \item being dlt is preserved under the minimal model program. 
\end{enumerate}
These properties also hold if we define ``dlt$=$lc$+$simple singularities".
\end{rem}

\subsection{Models} We define log birational models, weak lc models, log minimal models, and good minimal models of foliations similarly as the way they are defined for usual lc pairs.

\begin{defn}[Models]\label{defn: models}
Let $(X,\Ff,B)/U$ be an lc f-triple, $\phi: X\dashrightarrow X'$ a birational map over $U$, $E:=\Exc(\phi^{-1})$ the reduced $\phi^{-1}$-exceptional divisor, $\Ff':=\phi_*\Ff$, and $B':=\phi_*B+E^{\Ff'}$.
\begin{enumerate}
    \item $(X',\Ff',B')/U$ is called a \emph{log birational model} of $(X,\Ff,B)/U$. 
    \item $(X',\Ff',B')/U$ is called a \emph{weak lc model} of $(X,\Ff,B)/U$ if 
\begin{enumerate}
\item $(X',\Ff',B)/U$ is a log birational model of $(X,\Ff,B)/U$, 
    \item $K_{\Ff'}+B'$ is nef$/U$, and
    \item for any prime divisor $D$ on $X$ which is exceptional over $X'$, $a(D,\Ff,B)\leq a(D,\Ff',B')$.
\end{enumerate}
\item If $\rk\Ff=\dim X-1$ or $\dim X=3$ and $\rk\Ff=1$, then $(X',\Ff',B')/U$ is called a \emph{log minimal model} of $(X,\Ff,B)/U$ if
\begin{enumerate}
    \item $(X',\Ff',B')/U$ is a weak lc model of $(X,\Ff,B)/U$,
    \item $(X',\Ff',B')$ is $\Qq$-factorial dlt, and
    \item for any prime divisor $D$ on $X$ which is exceptional over $X'$, $a(D,\Ff,B)<a(D,\Ff',B')$.
\end{enumerate}
We remark that when $\dim X\leq 3$, by \cite[Theorem 11.3]{CS21} and Lemma \ref{lem: simple singularities cyc}, $X'$ is klt and $\Ff'$ has non-dicritical singularities.
\item If $\rk\Ff=\dim X-1$ or $\dim X=3$ and $\rk\Ff=1$, then $(X',\Ff',B')/U$ is called a \emph{good minimal model} of $(X,\Ff,B)/U$ if
\begin{enumerate}
        \item $(X',\Ff',B)/U$ is a log minimal model of $(X,\Ff,B)/U$, and
        \item $K_{\Ff'}+B'$ is semi-ample$/U$.
\end{enumerate}
\end{enumerate}
\end{defn}

\subsection{Minimal model program for foliations}

\begin{lem}[{cf. \cite[Proposition 8.1]{CS20}}]\label{lem: cs20 7.1 twist by nef}
Let $(X,\Ff,B)$ be an $f$-triple such that $\dim X=3$, $\rk\Ff=1$, and $(X,\Ff,B)$ is dlt. Let $H$ be a nef Cartier divisor on $X$ and $m\geq 2$ an integer such that $K_{\Ff}+B+mH$ is pseudo-effective. Then we may run a $(K_\Ff+B+mH)$-MMP 
$$X:=X_1\dashrightarrow X_2\dashrightarrow\dots\dashrightarrow X_i\dashrightarrow\dots,$$
and any such $(K_\Ff+B+mH)$-MMP satisfies the following. Let $\Ff_i$ be the induced foliation of $\Ff$ on $X_i$, and $B_i,H_i$ the images of $B,H$ on $X_i$ respectively. Then:
\begin{enumerate}
\item Each step of this $(K_\Ff+B+mH)$-MMP is $H$-trivial. In particular, $H_i$ is a nef Cartier divisor for each $i$.
\item $(X_i,\Ff_i,B_i)$ is dlt for each $i$.
\item The MMP terminates with a model $X_n$, such that $K_{\Ff_n}+B_n+mH_n$ is nef.
\end{enumerate}
We remark that $(X,\Ff,B+mH)$ may not be lc.
\end{lem}
\begin{proof}
First we prove (1)(2) by induction on $i$. Suppose we have already constructed a sequence of steps of a $(K_\Ff+B+mH)$-MMP 
$$X:=X_1\dashrightarrow X_2\dashrightarrow\dots\dashrightarrow X_k,$$
such that
\begin{itemize}
\item Each step of this partial $(K_\Ff+B+mH)$-MMP is  $H$-trivial. In particular, $H_i$ is a nef Cartier divisor for each $i\leq k$, and
\item $(X_i,\Ff_i,B_i)$ is dlt for each $i\leq k$.
\end{itemize}
If $K_{\Ff_k}+B_k+mH_k$ is nef then we are done. Otherwise, there exists a $(K_{\Ff_k}+B_k+mH_k)$-negative extremal ray $R$. Since $H_k$ is nef, any such $R$ is a $(K_{\Ff_k}+B_k)$-negative extremal ray. By \cite[Theorem 2.36]{CS20}, we may find an $\Ff$-invariant curve $C$ spanning $R$ such that $-2\leq (K_{\Ff_k}+B_k)\cdot C<0$. Since 
$(K_{\Ff_k}+B_k+mH_k)<0$, $m\geq 2$, and $H_k$ is nef Cartier, $H_k\cdot C=0$.

Since $(X_k,\Ff_k,B_k)$ is lc, no component of $B_k$ is $\Ff_k$-invariant. Thus $B_k\cdot C\geq 0$, so $K_{\Ff_k}\cdot C<0$. By \cite[Theorems 6.2, 6.3]{CS20}, the contraction associated to $R$ exists and, if the contraction is small, the flip exists. We let this step of the MMP be $X_k\dashrightarrow X_{k+1}$, $\Ff_{k+1}$ the induced foliation on $X_{k+1}$, and let $B_{k+1},H_{k+1}$ the images of $B_k,H_k$ on $X_{k+1}$ respectively. By \cite[Theorems 6.2, 6.3]{CS20}, $X_{k+1}$ has simple singularities. By \cite[Lemma 2.10]{CS20}, $(X_{k+1},\Ff_{k+1},B_{k+1})$ is lc, so $(X_{k+1},\Ff_{k+1},B_{k+1})$ is dlt. Since $H_k\cdot C=0$, $X_k\dashrightarrow X_{k+1}$ is $H_k$-trivial. In particular, $H_{k+1}$ is a nef Cartier divisor. (1)(2) follows by induction on $k$.

Now we prove (3). By (1), 
$$X_1\dashrightarrow X_2\dashrightarrow\dots\dashrightarrow $$
is a sequence of steps of a $(K_{\Ff}+B)$-MMP as well. By \cite[Theorem 7.1]{CS20}, this MMP terminates with a model $X_n$. Since  $K_{\Ff}+B+mH$ is pseudo-effective, $K_{\Ff_n}+B_n+mH_n$ is pseudo-effective, so $K_{\Ff_n}+B_n+mH_n$ is nef. This implies (3).
\end{proof}

\begin{lem}\label{lem: rel version 8.1 CS20}
Let $(X,\Ff,B)/U$ be an f-triple such that $\dim X=3$, $\rk\Ff=1$, $U$ is projective, $(X,\Ff,B)$ is dlt, and $K_{\Ff}+B$ is pseudo-effective$/U$. Then we may run a $(K_{\Ff}+B)$-MMP$/U$ and any such MMP$/U$ terminates with a log minimal model $(X',\Ff',B')/U$ of $(X,\Ff,B)/U$.
\end{lem}
\begin{proof}
Since $U$ is projective, $X$ is projective. Let $H$ be a general sufficiently ample divisor on $U$ and $\pi: X\rightarrow U$ the induced projective morphism. Then by Lemma \ref{lem: cs20 7.1 twist by nef}, we may run a $(K_{\Ff}+B+\pi^*H)$-MMP and any such MMP terminates with a model $X'$ satisfying the following. Let $\Ff'$ be the induced foliation of $\Ff$ on $X'$ and let $B',A'$ be the images of $B,\pi^*H$ on $X'$ respectively. Then
\begin{itemize}
\item Each step of this $(K_\Ff+B+\pi^*H)$-MMP is $\pi^*H$-trivial. In particular, $A'$ is a nef Cartier divisor,
\item $(X',\Ff',B')$ is dlt, and
\item $K_{\Ff'}+B'+A'$ is nef.
\end{itemize}
 Since $H$ is sufficiently ample, by the length of extremal rays \cite[Theorem 2.36]{CS20}, this MMP is an MMP$/U$, hence there exists an induced projective morphism $\pi': X'\rightarrow U$ and $A'=\pi'^*H$. Thus $(X',\Ff',B')/U$ is a log minimal model of $(X,\Ff,B)/U$.
\end{proof}

\begin{thm}\label{thm: fdlt modification}
Let $(X,\Ff,B)$ be an lc f-triple such that one of the following holds:
\begin{enumerate}
\item $\dim X\leq 3$ and $\rk\Ff=\dim X-1$.
\item $\dim X=3,\rk\Ff=1$, and $X$ is projective.
\end{enumerate}
Then $(X,\Ff,B)$ has a dlt model.
\end{thm}
\begin{proof}
(1) follows from \cite[Theorem 2.4]{SS22} (see \cite[Theorem 8.1]{CS21} for the projective case). Thus we may assume that $\dim X=3, \rk\Ff=1$, and $X$ is projective. By Theorem \ref{thm: existence of resolution}, there exists a foliated log resolution $h: W\rightarrow X$ of $(X,\Ff,B)$. Let $\Ff_W:=h^{-1}\Ff$, $E$ the reduced exceptional divisor of $h$, and $B_W:=h^{-1}_*B+E^{\Ff_W}$. Then $\Ff_W$ has simple singularities. By Lemma \ref{lem: corank log smooth imply flc}, $(W,\Ff_W,B_W)$ is lc. By construction, we may write 
$$K_{\Ff_W}+B_W=h^*(K_{\Ff}+B)+F$$
such that $F\geq 0$. The for any $h$-exceptional prime divisor $D$, $D$ is a comopnent of $\Supp F$ if and only if  $a(D,\Ff,B)>-\epsilon_{\Ff}(D)$.

We run a $(K_{\Ff_W}+B_W)$-MMP$/X$. By Lemma \ref{lem: rel version 8.1 CS20}, this MMP terminates with a log minimal model $(Y,\Ff_Y,B_Y)/X$ of $(W,\Ff_W,B_W)/X$. Since $K_{\Ff_W}+B_W\sim_{\mathbb R,X}F\geq 0$, by the negativity lemma, $K_{\Ff_Y}+B_Y\sim_{\mathbb R,X}0$. Thus we have $K_{\Ff_Y}+B_Y=f^*(K_{\Ff}+B)$, where $f: Y\rightarrow X$ is the induced birational morphism. By construction, $(Y,\Ff_Y,B_Y)$ is a dlt model of $(X,\Ff,B)$.
\end{proof}

The following theorem is a summary of the results in \cite{CS20,CS21,SS22}.

\begin{thm}\label{thm: can run mmp}
Let $(X,\Ff,B)/U$ be a $\Qq$-factorial lc f-triple such that $X$ is klt. Then we can run a $(K_{\Ff}+B)$-MMP$/U$ and any sequence of steps of a $(K_{\Ff}+B)$-MMP$/U$ will terminate with either a Mori fiber space$/U$ or a weak lc model of $(X,\Ff,B)/U$ if one of the following holds:
\begin{itemize}
    \item $\dim X\leq 3$, $\rk\Ff=\dim X-1$, and $(X,\Ff,B)$ is dlt.
    \item $\dim X=3$, $\rk\Ff=1$, and $U$ is projective.
\end{itemize}
Moreover, for any f-triple $(X',\Ff',B')/U$ constructed in a sequence of steps of a $(K_{\Ff}+B)$-MMP$/U$,
\begin{enumerate}
    \item $X'$ is $\Qq$-factorial klt, and
    \item if $(X,\Ff,B)$ is dlt, then $(X',\Ff',B')$ is dlt. 
\end{enumerate}
In particular, any sequence of steps of a $(K_{\Ff}+B)$-MMP$/U$ will terminate with either a Mori fiber space$/U$ or a log minimal model of $(X,\Ff,B)/U$.
\end{thm}
\begin{proof}
When $\dim X\leq 3$, $\rk\Ff=\dim X-1$, and $(X,\Ff,B)$ is dlt, the theorem follows from \cite[Theorem 0.5 and Corollary 2.3]{SS22} and \cite[Theorem 11.3]{CS21}, so we only need to prove the case when  $\dim X=3$ and $\rk\Ff=1$. In this case, the proof is similar to \cite[Proof of Theorem 8.8]{CS20}. We give a full proof here for the reader's convenience.

Since $(X,\Ff,B)$ is lc, no component of $B$ is $\Ff$-invariant.  By \cite[Theorem 2.36]{CS20}, there are countably many $(K_{\Ff}+B)$-negative extremal rays and all of them are spanned by $\Ff$-invariant curves.  

Let $H$ be a general sufficiently ample divisor on $U$ and $\pi: X\rightarrow U$ the induced projective morphism. If all $(K_{\Ff}+B)$-negative extremal rays are not $(K_{\Ff}+B)$-negative extremal rays$/U$, then by the length of extremal rays \cite[Theorem 2.36]{CS20}, $K_{\Ff}+B+\pi^*H$ is nef and we are done. Thus we may assume that there exists a $(K_{\Ff}+B)$-negative extremal ray$/U$ $R$ which is spanned by an $\Ff$-invariant curve $C$.

Since no components of $B$ are $\Ff$-invariant, $R$ is also a $K_{\Ff}$-negative extremal ray$/U$.  By \cite[Theorem 7.1]{CS20}, any sequence of $K_{\Ff}$-flips terminates, hence we only need to show that we can run a step of a $(K_{\Ff}+B)$-MMP$/U$ $(X,\Ff,B)\dashrightarrow (X',\Ff',B')$, and (1) and (2) hold. When $(X,\Ff,B)$ is dlt, the theorem follows from \cite[Theorems 6.2 and 6.3]{CS20}, so we may assume that $(X,\Ff,B)$ is not dlt. We only need to show that there exists a contraction of $R$ and the associated divisorial contraction, flip, or Mori fiber space exists, and (1) and (2) hold.

If $\loc R=X$, then the supporting function of $R$ is semi-ample by  \cite[Lemma 9.1]{CS20}, hence defines a contraction, so we get a Mori fiber space and we are done. Thus we may assume that $\loc R\not=X$. Then $\dim\loc R=1$ or $2$. By \cite[Corollary 8.4]{CS20}, $\Ff$ has canonical singularities in a neighborhood of $C$. If $\dim\loc R=1$, then by \cite[Theorem 8.7(1)]{CS20} and Lemma \ref{lem: simple singularities cyc}, the flip $X\dashrightarrow X^+$ associated to $R$ exists and $X^+$ is $\Qq$-factorial klt. Here we remark that \cite[Theorem 8.7(1)]{CS20} only mentions that $X^+$ is klt, but $X^+$ is also $\Qq$-factorial. In fact, $X^+$ is constructed in the following way:
\begin{center}$\xymatrix{
Y\ar@{-->}[r]^p\ar@{->}[d] & Y^+\ar@{-->}[d]\\
X & X^+.
}$
\end{center}
where $Y$ is $\Qq$-factorial and $\Gg:=p^{-1}\Ff$ has simple singularities \cite[Page 68, Lines 24-25]{CS20}, $Y\dashrightarrow Y^+$ is a $K_{\Gg}$-MMP \cite[Page 68, Line -9]{CS20} so  $Y^+$ is $\Qq$-factorial klt by Lemma \ref{lem: rel version 8.1 CS20}, and $Y^+\dashrightarrow X+$ is an MMP of an lc pair on $Y^+$ \cite[Page 69, Line 20]{CS20}, so $X^+$ is $\Qq$-factorial.

Thus we may assume that $\dim\loc R=2$. Then $D=\loc R$ is a divisor. By \cite[Theorem 6.2]{CS20}, the divisorial contraction $X\rightarrow Y$ associated to $R$ exists, and by \cite[Lemma 6.1; Page 60, Line -2; Page 61, Line 12-14]{CS20}, $X\rightarrow Y$ is a $(K_X+\lambda D)$-MMP where $\lambda$ is the lc threshold of $D$ with respect to $X$. Since all lc centers of $(X,\lambda D)$ are contained in $D$, $Y$ is $\Qq$-factorial klt.
\end{proof}

\begin{proof}[Proof of Theorem \ref{thm: abundance num dim 0}]
By Theorem \ref{thm: fdlt modification}, we may assume that $(X,\Ff,B)$ is dlt, $X$ is $\Qq$-factorial klt, and $\Ff$ is non-dicritical. By Theorem \ref{thm: can run mmp}, we can run a $(K_{\Ff}+B)$-MMP which terminates with a log minimal model $(X',\Ff',B')$ of $(X,\Ff,B)$ such that $\kappa_{\sigma}(K_{\Ff'}+B')=0$. By \cite[Chapter V, Proposition 2.7(6)]{Nak04}, $K_{\Ff'}+B'\equiv 0$. By Theorem \ref{thm: r coefficient nt abundance intro}, $K_{\Ff'}+B'\sim_{\mathbb R}0$. Thus $(X',\Ff',B')$ is a good minimal model of $(X,\Ff,B)$.  In particular, $\kappa_{\iota}(K_{\Ff}+B)=\kappa_{\iota}(K_{\Ff'}+B')=0$.
\end{proof}

Finally, we shall prove the existence of good minimal models for foliated pairs whose log canonical divisor is very exceptional, which is an analogue of \cite[Theorem 1.8]{Bir12}. This result will not be used in the rest of this paper, but it could be useful for future applications.

\begin{defn}[{cf. \cite[Definition 3.2]{Sho03}, \cite[Definition 3.1]{Bir12}}]
Let $f: X\rightarrow Y$ be a contraction between normal varieties, $D$ an $\Rr$-divisor on $X$, and $V\subset X$ a closed subset. We say that $V$ is \emph{vertical}/$Y$  if $f(V)$ is a proper subset of $Y$. We say that $V$ is \emph{horizontal}/$Y$  if $f(V)=Y$. We say that $D$ is \emph{very exceptional}/$Y$ if $D$ is vertical/$Y$, and for any prime divisor $P$ on $Y$, there exists a prime divisor $Q$ on $X$ which is not a component of $D$ but $f(Q)=P$, i.e. over the generic point of $P$ we have 
$$\Supp f^*P\not\subset\Supp D.$$
\end{defn}

\begin{thm}\label{thm: mmp very exceptional}
Let $(X,\Ff,B)/U$ be a $\Qq$-factorial lc f-triple such that $X$ is klt. Assume that there exists an $\Rr$-divisor $E\geq 0$ on $X$ such that $K_{\Ff}+B\sim_{\mathbb R,U}E$ and $E$ is very exceptional over $U$. Assume that one of the following holds:
\begin{itemize}
    \item $\dim X\leq 3$, $\rk\Ff=\dim X-1$, and $(X,\Ff,B)$ is dlt.
    \item $\dim X=3$, $\rk\Ff=1$, and $U$ is projective.
\end{itemize}
Then we can run a $(K_{\Ff}+B)$-MMP$/U$ and any sequence of steps of a $(K_{\Ff}+B)$-MMP$/U$ will terminate with a weak lc model $(X',\Ff',B')/U$ of $(X,\Ff,B)/U$ such that $K_{\Ff'}+B'\sim_{\mathbb R,U}0$. Moreover, 
\begin{enumerate}
    \item the image of $E$ on $X'$ is $0$, and
    \item if $(X,\Ff,B)$ is dlt, then $(X',\Ff',B')/U$ is a good minimal model of $(X,\Ff,B)/U$.
\end{enumerate}
\end{thm}
\begin{proof}
Since $\kappa_{\iota}(X/U,K_{\Ff}+B)=\kappa_{\iota}(X/U,E)\geq 0$, by Theorem \ref{thm: can run mmp},  we can run a $(K_{\Ff}+B)$-MMP$/U$ and any sequence of steps of a $(K_{\Ff}+B)$-MMP$/U$ will terminate with a weak lc model $(X',\Ff',B')/U$ of $(X,\Ff,B)/U$, and  $(X',\Ff',B')/U$ is a log minimal model of $(X,\Ff,B)/U$ if $(X,\Ff,B)$ is dlt. Let $E'$ be the image of $E$ on $X'$, then $K_{\Ff'}+B'\sim_{\mathbb R,U}E'\geq 0$ is nef over $U$ and very exceptional over $U$. By \cite[Lemma 3.3]{Bir12}, $E'=0$ and we are done.
\end{proof}

\section{Canonical bundle formulas for foliations}\label{sec: Canonical bundle formula}

In this section we establish the canonical bundle formula for foliations (Theorem \ref{thm: cbf ftriple}). Before we prove the main theorem, we introduce some basic definitions and notations.

\subsection{Lc-trivial fibration for foliations}

\begin{defn}[Discrepancy $\bb$-divisors]
Let $(X,\Ff,B,\Mm)$ be a sub-gfq. We define $\bb$-divisors
$\Aa(X,\Ff,B,\Mm)$ and $\Aa^*(X,\Ff,B,\Mm)$ in the following way. For any birational morphism $f: Y\rightarrow X$, we define
$$\Aa(X,\Ff,B,\Mm)_Y:=K_{\Ff_Y}+\Mm_Y-f^*(K_\Ff+B+\Mm_X) \text{ and }\Aa^*(X,\Ff,B,\Mm)_Y:=\Aa(X,\Ff,B,\Mm)^{>-1}_Y,$$
where $\Ff_Y:=f^{-1}\Ff$.
\end{defn}

\begin{defn}\label{defn: lc trivial fibration foliation}
Let $(X,\Ff,B)/U$ be an f-sub-triple and $f: X\rightarrow Z$ a contraction$/U$. We say that $f: (X,\Ff,B,\Mm)\rightarrow Z$ is an \emph{lc-trivial fibration$/U$} if
\begin{enumerate}
\item there exists a foliation $\Ff_Z$ on $Z$ such that $\Ff=f^{-1}\Ff_Z$,
\item $(X,\Ff,B,\Mm)$ is sub-lc over the generic point of $Z$,
\item $\rk f_*\mathcal{O}_X(\lceil\Aa^*(X,\Ff,B,\Mm)\rceil)=1$, and
    \item $K_{\Ff}+B+\Mm_X\sim_{\mathbb R,Z}0$.
\end{enumerate}
\end{defn}

\subsection{Equidimensional reduction}

We need the following definitions and results on equidimensional reduction for generalized pairs and toroidal generalized pairs.

\begin{defn}[{cf. \cite[Definition 2.1]{ACSS21}}]\label{defn: toroidal g-pairs}
Let $(X,\Sigma_X,\Mm)/U$ be a g-pair. We say that $(X,\Sigma_X,\Mm)$ is \emph{toroidal} if $\Sigma_X$ is a reduced divisor, $\Mm$ descends to $X$, and for any closed point $x\in X$, there exists a toric variety $X_{\sigma}$, a closed point $t\in X_{\sigma}$, and an isomorphism of complete local algebras 
$$\phi_x:\widehat{\mathcal{O}}_{X,x}\cong\widehat{\mathcal{O}}_{X_\sigma,t}.$$
such that the ideal of $\Sigma_X$ maps to the invariant ideal of $X_{\sigma}\backslash T_{\sigma}$, where $T_\sigma\subset X_\sigma$ is the maximal torus of $X_{\sigma}$. Any such pair $(X_\sigma, t)$ will be called as a \emph{local model} of $(X,\Sigma_X,\Mm)$ at $x\in X$.

Let $(X,\Sigma_X,\Mm)/U$ and $(Z,\Sigma_Z,\Mm^Z)/U$ be toroidal g-pairs and $f: X\rightarrow Z$ a surjective morphism$/U$. We say that $f: (X,\Sigma,\Mm)\rightarrow (Z,\Sigma,\Mm^Z)$ is a \emph{toroidal} morphism/is \emph{toroidal}, if for every closed point $x\in X$, there exist a local model $(X_\sigma,t)$ of $(X,\Sigma_X,\Mm)$ at $x$, a local model $(Z_{\tau},s)$ of $(Z,\Sigma_Z,\Mm^Z)$ at $z:=f(x)$, and a toric morphism $g: X_\sigma\to Z_{\tau}$, so that the diagram of algebras
\begin{center}$\xymatrix{
    \widehat{\mathcal{O}}_{X,x}\ar@{->}[r]^{\cong}  &     \widehat{\mathcal{O}}_{X_{\sigma},t} \\
     \widehat{\mathcal{O}}_{Z,z}\ar@{->}[r]^{\cong}\ar@{->}[u] & \widehat{\mathcal{O}}_{Z_{\tau},s}\ar@{->}[u]
}$
\end{center}
commutes. Here the vertical maps are the algebra homomorphisms induced by $f$ and $g$ respectively. 
\end{defn}

\begin{defn}[Property $(*)$, cf. {\cite[Definition 2.13]{ACSS21}}]
 Let $(X,B)$ be an lc pair and $f: X\rightarrow Z$ a contraction. We say that $(X,B)/Z$ satisfies \emph{Property $(*)$} if
 \begin{enumerate}
     \item there exists a reduced divisor $\Sigma_Z$ on $Z$ such that $(Z,\Sigma_Z)$ is log smooth,
     \item the vertical part of $B$ coincides with $f^{-1}(\Sigma_Z)$, and
     \item for any closed point $z\in Z$ and for any reduced divisor $\Sigma\geq\Sigma_Z$ such that $(Z,\Sigma)$ is log smooth near $z$, we have that $(X,B+f^*(\Sigma-\Sigma_Z))$ is lc near $f^{-1}(z)$.
 \end{enumerate}
\end{defn}

\begin{defthm}\label{defthm: weak ss reduction}
Let $X$ be a normal quasi-projective variety, $X\rightarrow U$ a projective morphism, $X\rightarrow Z$ a contraction$/U$, $B$ an $\Rr$-divisor on $X$, $\Mm$ an NQC$/U$ $\bb$-divisor on $X$, $D$ a prime divisor over $X$, and $D_Z$ a prime divisor over $Z$. Then there exist a toroidal g-pair  $(X',\Sigma_{X'},\Mm)/U$, a log smooth pair $(Z',\Sigma_{Z'})$, and a commutative diagram
 \begin{center}$\xymatrix{
X'\ar@{->}[r]^{h}\ar@{->}[d]_{f'}& X\ar@{->}[d]^{f}\\
Z'\ar@{->}[r]^{h_Z} & Z\\
}$
\end{center}
satisfying the following.
\begin{enumerate}
\item $h$ and $h_Z$ are projective birational morphisms.
\item $f': (X',\Sigma_{X'},\Mm)\rightarrow (Z',\Sigma_{Z'})$ is a toroidal morphism. 
\item $\Supp(h^{-1}_*B)\cup\Supp\Exc(h)$ is contained in $\Supp\Sigma_{X'}$,
\item $X'$ has at most toric quotient singularities. 
\item $f'$ is equidimensional.
\item $\Mm$ descends to $X'$.
\item $X'$ is $\Qq$-factorial klt.
\item The center of $D$  on $X'$ and the center of $D_Z$ on $Z'$ are divisors.
\end{enumerate}
We call any such $f': (X',\Sigma_{X'},\Mm)\rightarrow (Z',\Sigma_{Z'})$ (together with $h$ and $h_Z$) which satisfies (1-7) to be an \emph{equidimensional model}$/U$ of $f: (X,B,\Mm)\rightarrow Z$. 
\end{defthm}
\begin{proof}
Possibly replacing $Z$ with a high model and replacing $X$ accordingly, we may assume that $D_Z$ is a divisor on $Z$, hence the center of $D_Z$ on any high model of $Z$ is a divisor.  Possibly replacing $X$ with a high model, we may assume that $\Mm$ descends to $X$ and $D$ is a divisor on $X$. Therefore, $\Mm$ descends to any high model of $X$ and the center of $D$ on any high model of $X$ is a divisor. Now the theorem follows from \cite[Theorem 2.2]{ACSS21}, which in turn follows from \cite[Theorem 2.1 and Proposition 4.4]{AK00}. We also refer the reader to \cite[Theorem B.6]{Hu20} for a more detailed explanation.
\end{proof}

\begin{lem}\label{lem: toroidal lct=1}
    Let $f: (X,\Sigma_X)\rightarrow (Z,\Sigma_Z)$ be a toroidal morphism such that $(Z,\Sigma_Z)$ is log smooth. Let $(X,B)$ be a sub-pair such that $\Sigma_X\geq B$
    and $B^{\geq 0}$ is horizontal$/Z$. 

    Then for any prime divisor $D$ on $X$, $(X,B+f^{-1}(D))$ is lc over the generic point of $D$.
\end{lem}
\begin{proof}
If $D\subset\Sigma_Z$, then $\Sigma_X\geq f^{-1}(D)$ and all components of $f^{-1}(D)$ are vertical$/Z$. Since $B^{\geq 0}$ is horizontal$/Z$, $\Sigma_X\geq B+f^{-1}(D)$. Thus $(X,B+f^{-1}(D))$ is lc. Thus we may assume that $D\not\subset\Sigma_Z$.

By \cite[Proposition 2.16]{ACSS21}, $(X,\Sigma_X)$ satisfies Property $(*)$. By \cite[Lemma 2.14(2)]{ACSS21}, $(X,\Sigma_X+f^*D)$ is lc over the generic point of $D$. Since $f^*D\geq f^{-1}(D)$ over the generic point of $D$ and $\Sigma_X\geq B$, $(X,B+f^{-1}(D))$ is lc over the generic point of $D$.
\end{proof}

\subsection{Bounded families}

We recall the definition of bounded families.

\begin{defn}[Bounded families]\label{defn: bdd}
A collection of varieties $\mathcal{D}$ is
said to be a \emph{bounded family} if there exists a projective morphism $h\colon \mathcal{Z}\rightarrow S$
of schemes of finite type such that
each $X\in \mathcal{D}$ is isomorphic to $\mathcal{Z}_s$ 
for some closed point $s\in S$.
\end{defn}

\subsection{Canonical bundle formulas}

\begin{setup}\label{set-up: cbf setup 1}
$X,\Ff,B,\Mm,U,f,Z,\Ff_Z$ satisfy the following:
\begin{enumerate}
    \item $(X,\Ff,B,\Mm)/U$ is a sub-gfq.
    \item $f: (X,\Ff,B,\Mm)\rightarrow Z$ is an lc-trivial fibration$/U$.
    \item $\Ff_Z$ is a foliation on $Z$ such that $f^{-1}\Ff_Z=\Ff$.
\end{enumerate}
\end{setup}

\begin{setup}\label{set-up: cbf setup 2}
$X,\Ff,B,\Mm,U,f,Z,\Ff_Z,f',X',\Sigma_{X'},Z',\Sigma_{Z'},h,h_Z,\Ff',\Ff_{Z'},R'$ satisfy the following:
\begin{enumerate}
\item $X,\Ff,B,\Mm,U,f,Z,\Ff_Z$ satisfy the conditions of Set-up \ref{set-up: cbf setup 1}.
\item $f': (X',\Sigma_{X'},\Mm)\rightarrow (Z',\Sigma_{Z'})$ together with $h: X'\rightarrow X$ and $h_Z: Z'\rightarrow Z$ is an equidimensional model$/U$ of $f: (X,B,\Mm)\rightarrow Z$, whose existence is guaranteed by Definition-Theorem \ref{defthm: weak ss reduction}.
\item $\Ff':=h^{-1}\Ff$, $\Ff_{Z'}:=h_Z^{-1}\Ff_Z$.
\item $$K_{\Ff'}+B'+\Mm_{X'}:=h^*(K_X+B+\Mm_X).$$
\item 
    $$R':=\sum_{P\mid P\text{ is an }\Ff_{Z'}\text{-invariant prime divisor}}(f'^*P-f'^{-1}(P))$$
\end{enumerate}
\begin{center}$\xymatrix{
(X',\Sigma_{X'},\Mm)\ar@{->}[r]^h\ar@{->}[d]_{f'} & X\ar@{->}[d]^{f}\\
(Z',\Sigma_{Z'})\ar@{->}[r]^{h_Z} & Z}$
\end{center}
\end{setup}

\begin{lem}\label{lem: relation of f and x on equi model}
       Conditions as in Set-up \ref{set-up: cbf setup 2}. Then:
       \begin{enumerate}
           \item
                  $K_{\Ff'/\Ff_{Z'}}=K_{X'/Z'}-R'.$
            \item  $\rk f'_*\mathcal{O}_{X'}(\lceil\Aa^*(X',B'-R',\Mm)\rceil)=1$.
            \item $K_{X'}+B'-R'+\Mm_{X'}\sim_{\mathbb R,Z'}0$.
            \item $(X',B'-R',\Mm)$ is sub-lc over the generic point of $Z'$.
       \end{enumerate}
\end{lem}
\begin{proof}
    (1) follows from  \cite[2.9]{Dru17}. Since $K_{\Ff}+B+\Mm_{X}\sim_{\mathbb R,Z}0$, $$K_{\Ff'}+B'+\Mm_{X'}=h^*(K_{\Ff}+B+\Mm_X)\sim_{\mathbb R,Z}0.$$ 
    Since $\rk f_*\mathcal{O}_X(\lceil\Aa^*(X,\Ff,B,\Mm)\rceil)=1$, $\rk f_*\mathcal{O}_X(\lceil\Aa^*(X',\Ff',B',\Mm)\rceil)=1$. By (1),
    $$K_{X'}+B'-R'+\Mm_{X'}=(K_{\Ff'}+B'+\Mm_{X'})+f'^*(K_{Z'}-K_{\Ff_{Z'}}),$$
    which implies (2). Moreover, since
    $$(K_{\Ff'}+B'+\Mm_{X'})+f'^*(K_{Z'}-K_{\Ff_{Z'}})\sim_{\mathbb R,Z'}K_{\Ff'}+B'+\Mm_{X'},$$
    we get (3). Since $K_{\Ff'}=K_{X'}$ over the generic point of $Z'$, we get (4).
\end{proof}

\begin{deflem}\label{deflem: bd and td}
    Conditions as in Set-up \ref{set-up: cbf setup 1}. For any prime divisor $D$ over $Z$, we define $b_D(X,\Ff,B,\Mm;\Ff_Z)$ and $t_D(X,\Ff,B,\Mm;\Ff_Z)$ in the following way. 

    By Definition-Theorem \ref{defthm: weak ss reduction}, we may choose $f',X',\Sigma_{X'},Z',\Sigma_{Z'},h,h_Z,\Ff',\Ff_{Z'},R'$ so that the conditions of Set-up \ref{set-up: cbf setup 2} are satisfied and $D$ is a divisor on $Z'$. Now we define
    $$b_D(X,\Ff,B,\Mm;\Ff_Z):=1-\sup\{t\mid (X',B'-R'+tf'^*D,\Mm)\text{ is sub-lc over the generic point of } D\},$$
    and
    \begin{align*}
    t_D(X,\Ff,B,\Mm;\Ff_Z):=&\epsilon_{\Ff_{Z'}}(D)-\Bigg\{t\Biggm|
    \begin{array}{r@{}l}
       (X',\Ff',B'+tf'^*D,\Mm)\text{ is sub-lc}\\ 
       \text{over the generic pont of }D 
    \end{array}\Bigg\}.
    \end{align*}
    Then $b_D(X,\Ff,B,\Mm;\Ff_Z)$ and $ t_D(X,\Ff,B,\Mm;\Ff_Z)$ are well-defined, i.e. they are independent of the choices of the equidimensional model$/U$ of $f: (X,B,\Mm)\rightarrow Z$ and are independent of $U$.
\end{deflem}
\begin{proof}
     By Definition-Theorem \ref{defthm: weak ss reduction}, we may choose $f',X',\Sigma_{X'},Z',\Sigma_{Z'},h,h_Z,\Ff',\Ff_{Z'},R'$ so that the conditions of Set-up \ref{set-up: cbf setup 2} are satisfied. Since $b_D(X,\Ff,B,\Mm;\Ff_Z)$ and $ t_D(X,\Ff,B,\Mm;\Ff_Z)$ are invariants which measure singularities, they are independent of the choice of $U$.

We let $f'': (X'',\Sigma_{X''},\Mm)\rightarrow (Z'',\Sigma_{Z''})$ together with $g: X''\rightarrow X$ and $g_Z: Z''\rightarrow Z$ be another equidimensional model$/U$ of $f: (X,B,\Mm)\rightarrow Z$ such that $D$ is on $Z''$. Since the question is local near the generic point of $D$, possibly by shrinking $Z'$ and $Z''$ to a neighborhood of $D$, we may assume that $g_Z=h_Z$, and $Z'=Z''$. We let $\Ff'':=g^{-1}\Ff$,
 $$K_{\Ff''}+B''+\Mm_{X''}:=g^*(K_X+B+\Mm_X),$$
and
$$R'':=\sum_{P\mid P\text{ is an }\Ff_{Z'}\text{-invariant prime divisor}}(f''^*P-f''^{-1}(P)).$$
We let $p: W\rightarrow X'$ and $q: W\rightarrow X''$ be a common resolution and $\Ff_W:=p^{-1}\Ff'$. Then $\Ff_W=q^{-1}\Ff''$. Since $(X',\Ff',B')$ and $(X'',\Ff'',B'')$ are both crepant to $(X,\Ff,B)$, we have
$$K_W+B_W+\Mm_W:=p^*(K_{\Ff'}+B'+\Mm_{X'})=q^*(K_{\Ff''}+B''+\Mm_{X''}).$$
Therefore,
\begin{align*}
&\sup\{t\mid (X',\Ff',B'+tf'^*D,\Mm)\text{ is sub-lc over the generic point of } D\}\\
=&\sup\{t\mid (W,\Ff_W,B_W+t(f'\circ p)^*D,\Mm)\text{ is sub-lc over the generic point of } D\}\\
=&\sup\{t\mid (W,\Ff_W,B_W+t(f''\circ q)^*D,\Mm)\text{ is sub-lc over the generic point of } D\}\\
=&\sup\{t\mid (X'',\Ff'',B''+tf''^*D,\Mm)\text{ is sub-lc over the generic point of } D\}.
\end{align*}
So $t_D(X,\Ff,B,\Mm;\Ff_Z)$ is independent of the choices of the equidimensional model$/U$ of $f: (X,B,\Mm)\rightarrow Z$.

By Lemma \ref{lem: relation of f and x on equi model}(1), we have
\begin{align*}
    p^*(K_{X'}+B'-R'+\Mm_{X'})&=p^*(K_{\Ff'}+B'+\Mm_{X'}+f'^*(K_{Z'}-K_{\Ff_{Z'}}))\\
    &=q^*(K_{\Ff''}+B''+\Mm_{X''})+(f'\circ p)^*(K_{Z'}-K_{\Ff_{Z'}})\\
    &=q^*(K_{\Ff''}+B''+\Mm_{X''})+(f''\circ q)^*(K_{Z'}-K_{\Ff_{Z'}})\\
    &=q^*(K_{\Ff''}+B''+\Mm_{X''}+f''^*(K_{Z''}-K_{\Ff_{Z''}}))\\
    &=q^*(K_{X''}+B''-R''+\Mm_{X''}).
\end{align*}
Let $K_W+D_W+\Mm_W:=p^*(K_{X'}+B'-R'+\Mm_{X'})$, then
\begin{align*}
&\sup\{t\mid (X',B'-R'+tf'^*D,\Mm)\text{ is sub-lc over the generic point of } D\}\\
=&\sup\{t\mid (W,D_W+t(f'\circ p)^*D,\Mm)\text{ is sub-lc over the generic point of } D\}\\
=&\sup\{t\mid (W,D_W+t(f''\circ q)^*D,\Mm)\text{ is sub-lc over the generic point of } D\}\\
=&\sup\{t\mid (X'',B''-R''+tf''^*D,\Mm)\text{ is sub-lc over the generic point of } D\}.
\end{align*}
So $b_D(X,\Ff,B,\Mm;\Ff_Z)$ is independent of the choices of the equidimensional model$/U$ of $f: (X,B,\Mm)\rightarrow Z$.
\end{proof}

Now we are ready to establish the canonical bundle formula for foliations.

\begin{deflem}\label{deflem: discriminant and moduli part}
    Conditions as in Set-up \ref{set-up: cbf setup 1}. Further assume that $\Mm$ is semi-ample$/Z$.  By Definition-Theorem \ref{defthm: weak ss reduction}, we may choose $f',X',\Sigma_{X'},Z',\Sigma_{Z'},h,h_Z,\Ff',\Ff_{Z'},R'$ so that the conditions of Set-up \ref{set-up: cbf setup 2} are satisfied. We define two $\bb$-divisors $\Bb$ and $\Mm^Z$ on $Z$ in the following way. By Lemma \ref{lem: relation of f and x on equi model}, $f': (X',B'-R',\Mm)\rightarrow Z'$ is an lc-trivial fibration. By \cite[Theorem 2.23]{JLX22} (see \cite[Theorem 9]{Fil19} for the $\Qq$-coefficient case), there exist two $\bb$-divisors $\Bb$ and $\Mm^Z$ on $Z$ satisfying the following.
    \begin{enumerate}
        \item[(i)] $K_{X'}+B'-R'+\Mm_{X'}\sim_{\mathbb R}f'^*(K_{Z'}+\Bb_{Z'}+\Mm^Z_{Z'})$.
        \item[(ii)] $\Mm^Z$ is NQC$/U$.
        \item[(iii)] For any birational morphism $g_Z: Z''\rightarrow Z'$ and $g: X''\rightarrow X'$ such that the induced map $f'': X''\dashrightarrow Z''$ is a morphism, we let 
        $$K_{X''}+G''+\Mm_{X''}:=g^*(K_{X'}+B'-R'+\Mm_{X'}),$$ 
        then for any prime divisor $D$ on $Z''$, we have
        $$\mult_D\Bb_{Z''}=1-\sup\{t\mid (X'',G''+tf''^*D,\Mm)\text{ is sub-lc over the generic point of }D\}$$
        and
        $$K_{Z''}+\Bb_{Z''}+\Mm^Z_{Z''}=g_Z^*(K_{Z'}+\Bb_{Z'}+\Mm^Z_{Z'}).$$
        \item[(iv)] $\Bb$ is a uniquely determined $\Rr$-divisor, and the $\Rr$-linear equivalence class of $\Mm^Z$ is uniquely determined. 
    \end{enumerate}
     We call $\Bb$ as the \emph{discriminant part} of $f: (X,\Ff,B,\Mm)\rightarrow Z$ and call $\Mm^Z$ as the \emph{moduli part} of $f: (X,\Ff,B,\Mm)\rightarrow Z$. Then:
    \begin{enumerate}
      \item  $\Bb$ and $\Mm^Z$ are well-defined, i.e. $\Bb$ 
 and the $\Rr$-linear equivalence class of $\Mm^Z$ are independent of the choices of the equidimensional model$/U$ of $f: (X,B,\Mm)\rightarrow Z$.  
 \item $(Z,\Ff_Z,B_Z:=\Bb_Z,\Mm^Z)/U$ is a sub-gfq. We say that $(Z,\Ff_Z,B_Z,\Mm^Z)/U$ is a sub-gfq \emph{induced by a canonical bundle formula$/U$} of $f: (X,\Ff,B,\Mm)\rightarrow Z$.
    \end{enumerate}
\end{deflem}
\begin{proof}
By Definition-Lemma \ref{deflem: bd and td}, for any birational model $Z''$ of $Z$ and any prime divisor $D$ on $Z''$, $\mult_D\Bb_{Z''}=b_D(X,\Ff,B,\Mm;\Ff_Z)$ which is independent of the choice of the equidimensional model$/U$ of $f: (X,B,\Mm)\rightarrow Z$. Thus $\Bb$ is independent of the choice of the equidimensional model$/U$ of $f: (X,B,\Mm)\rightarrow Z$.

Since $K_{\Ff}+B+\Mm_X\sim_{\mathbb R,Z}0$, there exists an $\Rr$-divisor $L$ on $Z$ which is uniquely determined up to $\Rr$-linear equivalence, such that
$$K_{\Ff}+B+\Mm_X\sim_{\mathbb R}f^*L.$$
By condition (i), we have
$$K_{\Ff'}+B'+\Mm_{X'}\sim_{\mathbb R}f'^*(K_{\Ff_{Z'}}+\Bb_{Z'}+\Mm^Z_{Z'}).$$
Therefore, for any birational morphism $g_Z: Z''\rightarrow Z'$ with $\Ff_{Z''}:=g_Z^{-1}\Ff_{Z'}$, we have
$$\Mm^Z_{Z''}\sim_{\mathbb R}(h_Z\circ g_Z)^*L-K_{\Ff_{Z''}}-\Bb_{Z''}.$$
Thus $\Mm^Z_{Z''}$ is uniquely determined up to the choices of $L$ in its $\Rr$-linear equivalence class. Thus $\Mm^Z$ is uniquely determined up to $\Rr$-linear equivalence. This implies (1).

Moreover, we have
$$L=(h_Z)_*h_Z^*L\sim_{\mathbb R}(h_Z)_*(K_{\Ff_{Z'}}+\Bb_{Z'}+\Mm^Z_{Z'})=K_{\Ff_Z}+B_Z+\Mm^Z_Z,$$
so $K_{\Ff_Z}+B_Z+\Mm^Z_Z$ is $\Rr$-Cartier. By our condition (ii), $(Z,\Ff_Z,B_Z:=\Bb_Z,\Mm^Z)/U$ is a sub-gfq. This implies (2).
\end{proof}

\subsection{Proof of Theorem \ref{thm: cbf ftriple}}

In this subsection, we prove Theorem \ref{thm: cbf ftriple}. Indeed, we shall prove a more general version of Theorem \ref{thm: cbf ftriple} for gfqs with the extra condition that the nef part is $\bb$-semiample.

\begin{prop}\label{prop: cbf m is generic}
    Conditions as in Set-up \ref{set-up: cbf setup 1}. Further assume that $\Mm$ is semi-ample$/Z$. We let $\Bb$ be the discriminant part of $f: (X,\Ff,B,\Mm)\rightarrow Z$, $\Mm^Z$ the moduli part of $f: (X,\Ff,B,\Mm)\rightarrow Z$, and $(Z,\Ff_Z,B_Z:=\Bb_Z,\Mm^Z)/U$ a sub-gfq induced by a canonical bundle formula of $f: (X,\Ff,B,\Mm)\rightarrow Z$. 

    Then the class of $\Mm^Z$ only depends on $(X,B,\Mm)/U$ over the generic point of $Z$. More precisely, for any sub-gfq $(\bar Z,\Ff_{\bar Z},B_{\bar Z},\Mm^{\bar Z})/U$ induced by an lc-trivial fibration$/U$ $\bar f: (\bar X,\bar \Ff,\bar B,\bar \Mm)\rightarrow\bar Z$ of a sub-gfq $(\bar X,\bar\Ff,\bar B,\bar\Mm)/U$ such that $f=\bar f$, $B=\bar B$, and $\Mm=\bar\Mm$ over the generic point of $Z$, we have $\Mm^Z\sim_{\mathbb R,U}\Mm^{\bar Z}$.

In particular, for any integer $n$ such that $n(K_X+B+\Mm_X)\sim 0$ over the generic point of $Z$, we may choose $\Mm^Z$ such that $n(K_{\Ff}+B+\Mm_X)\sim nf^*(K_{\Ff_Z}+B_Z+\Mm^Z_Z)$.
\end{prop}
\begin{proof}
 By Definition-Theorem \ref{defthm: weak ss reduction}, we may choose $f',X',\Sigma_{X'},Z',\Sigma_{Z'},h,h_Z,\Ff',\Ff_{Z'},R'$ so that the conditions of Set-up \ref{set-up: cbf setup 2} are satisfied.
 
     Since $R'$ is vertical$/Z'$, $\Mm^Z$ only depends on $(X',B',\Mm')$ over the generic point of $Z'$ (see \cite[Lemma 3.5]{Bir19} for the $\Qq$-coefficient and $\Mm=\bm{0}$ case; the general case follows from the same lines of the proof). Since $K_{X'}+B'+\Mm_{X'}=h^*(K_X+B+\Mm_X)$ over the generic point of $Z$, and $Z'$ is birational to $Z$, the $\Rr$-linear equivalence class of $\Mm^Z$ only depends on $(X,B,\Mm)$ over the generic point of $Z$. The last argument of the proposition follows from \cite[Definition 2.19]{FS23} and \cite[7.5 Construction]{PS09}. 
\end{proof}

\begin{prop}\label{prop: cbf is really cbf}
Conditions as in Set-up \ref{set-up: cbf setup 1}. Further assume that $\Mm$ is semi-ample$/Z$ and $(X,\Ff,B,\Mm)$ is sub-lc. We let $\Bb$ be the discriminant part of $f: (X,\Ff,B,\Mm)\rightarrow Z$, $\Mm^Z$ the moduli part of $f: (X,\Ff,B,\Mm)\rightarrow Z$, and $(Z,\Ff_Z,B_Z:=\Bb_Z,\Mm^Z)/U$ a sub-gfq induced by a canonical bundle formula of $f: (X,\Ff,B,\Mm)\rightarrow Z$. Then:
\begin{enumerate}
\item $(Z,\Ff_Z,B_Z,\Mm^Z)$ is sub-lc.
\item If $B\geq 0$, then 
\begin{enumerate}
    \item $(Z,\Ff_Z,B_Z,\Mm^Z)$ is lc, and
    \item for any component $D$ of $B_Z$, $\mult_DB_Z=t_D(X,\Ff,B,\Mm;\Ff_Z)$.
\end{enumerate}
\end{enumerate}
\end{prop}
\begin{proof}
By Definition-Lemma \ref{deflem: bd and td}, we only need to show that for any prime divisor $D$ over $X$, 
\begin{itemize}
   \item $b_D(X,\Ff,B,\Mm;\Ff_Z)\leq\epsilon_{\Ff_Z}(D)$, and
   \item if $B\geq 0$ and $D$ is not exceptional over $Z$, then
   \begin{itemize}
       \item[$\cdot$]
    $b_D(X,\Ff,B,\Mm;\Ff_Z)\geq 0$, and
       \item[$\cdot$] $b_D(X,\Ff,B,\Mm;\Ff_Z)=t_D(X,\Ff,B,\Mm;\Ff_Z)$.
   \end{itemize}
\end{itemize}
 By Definition-Theorem \ref{defthm: weak ss reduction}, we may choose $f',X',\Sigma_{X'},Z',\Sigma_{Z'},h,h_Z,\Ff',\Ff_{Z'},R'$ so that the conditions of Set-up \ref{set-up: cbf setup 2} are satisfied and $D$ is a divisor on $Z'$.

 If $D$ is not $\Ff_{Z'}$-invariant, then $f'^{-1}(D)$ does not intersect $R'$ over the generic point of $D$. Since $\Sigma_{X'}\geq\Supp B'\geq B'$ and $(X',\Sigma_{X'},\Mm)$ is lc, $(X',B',\Mm)$ is lc. Thus we have
\begin{align*}
   &b_D(X,\Ff,B,\Mm;\Ff_Z)\\
   =&1-\sup\{t\mid (X',B'-R'+tf'^*D,\Mm)\text{ is sub-lc over the generic point of } D\}\\
   =&1-\sup\{t\mid (X',B'+tf'^*D,\Mm)\text{ is sub-lc over the generic point of } D\}\\
   \leq &1=\epsilon_{\Ff_{Z'}}(D).
\end{align*}
Moreover, since $D$ is not $\Ff_{Z'}$-invariant, no component of $f'^{-1}(D)$ is $\Ff'$-invariant. Thus
\begin{align*}
    &\sup\{t\mid (X',B'+tf'^*D,\Mm)\text{ is sub-lc over the generic point of } D\}\\
    =&\sup\{t\mid (X',\Ff',B'+tf'^*D,\Mm)\text{ is sub-lc over the generic point of } D\}\\
    =&1-t_D(X,\Ff,B,\Mm;\Ff_Z).
\end{align*}
In particular, if $D$ is not exceptional over $Z$, then
\begin{align*}
    \mult_DB_Z&=b_D(X,\Ff,B,\Mm;\Ff_Z)=t_D(X,\Ff,B,\Mm;\Ff_Z)\\
  &=1-\sup\{t\mid (X,\Ff,B+tf^*D,\Mm)\text{ is sub-lc over the generic point of } D\},
\end{align*}
so if $B\geq 0$, then $\mult_DB_Z\geq 0$.

If $D$ is $\Ff_{Z'}$-invariant, then any component of $f'^{-1}(D)$ is $\Ff'$-invariant. Since $(X',\Ff',B',\Mm)$ is sub-lc, for any component $C$ of $f'^{-1}(D)$, $\mult_CB'\leq 0$. By Lemma \ref{lem: toroidal lct=1}, $(X',B'+f'^{-1}(D),\Mm)$ is sub-lc over the generic point of $D$. Since
$(X',B'-R'+f'^*D,\Mm)=(X',B'+f'^{-1}(D),\Mm)$
over the generic point of $D$,  
$$\sup\{t\mid (X',B'-R'+tf'^*D,\Mm)\text{ is sub-lc over the generic point of } D\}\geq 1.$$
Thus $b_D(X,\Ff,B,\Mm;\Ff_Z)\leq 0=\epsilon_{\Ff_{Z'}}(D)$. Moreover, if $D$ is not exceptional over $Z$, then $f'^{-1}(D)$ is not exceptional over $X$. Therefore, if $B\geq 0$, then
$$\sup\{t\mid (X',B'-R'+tf'^*D,\Mm)\text{ is sub-lc over the generic point of } D\}=1,$$ 
and we have $$\mult_DB_Z=b_D(X,\Ff,B,\Mm;\Ff_Z)=t_D(X,\Ff,B,\Mm;\Ff_Z)=0.$$
The proposition follows.
\end{proof}

\begin{proof}[Proof of Theorem \ref{thm: cbf ftriple}]
This is a special case of Propositions \ref{prop: cbf m is generic} and \ref{prop: cbf is really cbf}.
\end{proof}

\subsection{Low-dimensional canonical bundle formulas}

In low dimensions, we have a more detailed understanding on the canonical bundle formula of foliations as Prokhorov-Shokurov's $\bb$-semiampleness conjecture holds. We have the following proposition which will be used later in this paper.

\begin{prop}\label{prop: low-dimension cbf}
Let $\Ii_0\subset [0,1]\cap\mathbb Q$ be a finite set. Then there exists a positive integer $I$ depending only on $\Ii_0$ satisfying the following. 

Let $(X,\Ff,B)$ be a projective lc f-triple, $f: X\rightarrow Z$ a contraction$/U$, and $\Ff_Z$ a foliation on $Z$ such that $\Ff=f^{-1}\Ff_Z$. We let $\Bb$ be the discriminant part of $f: (X,\Ff,B)\rightarrow Z$, $\Mm^Z$ the moduli part of $f: (X,\Ff,B)\rightarrow Z$, and $(Z,\Ff_Z,B_Z:=\Bb_Z,\Mm^Z)$ an lc gfq induced by a canonical bundle formula of $f: (X,\Ff,B)\rightarrow Z$ (whose existence is guaranteed by Proposition \ref{prop: cbf is really cbf}). Assume that
\begin{enumerate}
    \item $\dim X\leq 3$,
    \item $B^h\in\Ii_0$ where $B^h$ is the horizontal$/Z$ part of $B$, and
    \item one of the following holds:
    \begin{enumerate}
            \item $\dim Z=\dim X-1$.
    \item $\dim X=3,\dim Z=1$, and $B^h\not=0$.
    \item $\dim X=3,\dim Z=1$, and $\Mm^Z_Z\equiv 0$.
    \end{enumerate}
\end{enumerate}
Then we may choose $\Mm^Z$ such that $I\Mm^Z$ is base-point-free, and for any integer $n$ divisible by $I$ such that $n(K_X+B+\Mm_X)\sim 0$ over the generic point of $Z$, we have
$$n(K_{\Ff}+B+\Mm_X)\sim nf^*(K_{\Ff_Z}+B_Z+\Mm^Z_Z).$$
\end{prop}

\begin{proof}
Possibly by replacing $\Ii_0$ with $\Ii_0\cup\{1\}$ and replacing $(X,\Ff,B)$ with a dlt modification of $(X,\Ff,B)$ (whose existence is guaranteed by Theorem \ref{thm: fdlt modification}), we may assume that $(X,\Ff,B)$ is $\Qq$-factorial dlt.

If $\dim Z=\dim X-1$, then the general fiber of $f$ is either a rational curve or an elliptic curve. In this case, the proposition follows from Proposition \ref{prop: cbf m is generic}, \cite[Theorem 8.1]{PS09} and \cite{Kod63,Uen73}. If $\dim X=3,\dim Z=1$, and $B^h\not=0$, then the proposition follows from Proposition \ref{prop: cbf m is generic} and \cite[Theorem 1.4]{ABBDILW23}.

Therefore, we may assume that $\dim X=3,\dim Z=1,B^h=0$, and $\Mm^Z_Z\equiv 0$. In particular, $(X,B)$ is klt over the generic point of $Z$. Let $F$ be a general fiber of $f$. Then $B|_F=0$ and $$K_F=(K_X+B)|_F=(K_{\Ff}+B)|_F\sim_{\mathbb Q}0$$
since $K_{\Ff}=K_X$ over the generic point of $Z$. By the global ACC \cite[Theorem 1.5]{HMX14}, $F$ is $\epsilon$-lc for some $\epsilon>0$ depending only on $\Ii_0$.

We let $b$ be the minimal positive integer such that $bK_F\sim 0$, $\tilde F$ be a smooth model of the $\mathbb Z/(b)$-cover of $F$, and $\beta$ the second Betti number of $\tilde F$. 

If $F$ is not canonical, then by \cite[0.4]{Ale94}, $F$ belongs to a bounded family. Thus $\tilde F$ belongs to a bounded family, so $b$ and $\beta$ are bounded. Thus the proposition follows from Proposition \ref{prop: cbf m is generic} and \cite[Theorem 1.3]{Flo14} in this case.

If $F$ is canonical, then we let $F'\rightarrow F$ be the minimal resolution of $F$. By the classification of smooth Calabi-Yau surfaces, either $4K_{F'}\sim 0$ or $6K_{F'}\sim 0$. Thus $4K_{F}\sim 0$ or $6K_{F}\sim 0$, so $b\leq 6$. Moreover, by our construction, $K_{\tilde F}\sim 0$ and $\tilde F$ is smooth, so $\tilde F$ is either a K3 surface with $\beta=22$, or an abelian surface with $\beta=6$. Thus $\beta\leq 22$. Now the proposition follows from Proposition \ref{prop: cbf m is generic} and \cite[Theorem 1.3]{Flo14} in this case as well, and the proposition follows.
\end{proof}

\section{Index theorem for surfaces}\label{sec: Surface index theorem, global version}

\begin{lem}\label{lem: surface canonical 12 index}
Let $\Ff$ be a foliation of rank $1$ on a normal projective surface $X$ such that $\Ff$ is canonical and $K_{\Ff}\equiv 0$. Then $IK_{\Ff}\sim 0$ for some positive integer $I\leq 12$.
\end{lem}
\begin{proof}
By Theorem \ref{thm: r coefficient nt abundance intro}, $K_{\Ff}\sim_{\mathbb Q}0$. Let $f: Y\rightarrow X$ be a resolution of $\Ff$ and $\Gg$ the induced foliation of $\Ff$ on $Y$. Since $\Ff$ is canonical, $K_{\Gg}\geq f^*K_{\Ff}$, $\kappa(K_{\Gg})=\kappa(K_{\Ff})=0$. By \cite[Theorem 1]{Per05}, $|IK_{\Gg}|\not=\emptyset$ for some positive integer $I\leq 12$, hence $|IK_{\Ff}|\not=\emptyset$ and $IK_{\Ff}\sim 0$.
\end{proof}

\begin{lem}\label{lem: surface gloabl index with boundary}
Let $\Ii\subset [0,1]$ be a DCC set and $p$ a positive integer. Then there exists a finite set $\Ii_0$ and a positive integer $I$ depending only on $\Ii$ and $p$ satisfying the following. Suppose that $(X,\Ff,B,\Mm)$ is a projective lc gfq, such that $\dim X=2$, $\rk\Ff=1$, $B\in\Ii$, $p\Mm$ is base-point-free,  $K_{\Ff}+B+\Mm_X\sim_{\mathbb Q}0$, $B\not=0$ or $\Mm_X\not\equiv 0$, and  $X$ is klt or $\Mm=\bm{0}$, then:
\begin{enumerate}
    \item $B\in\Ii_0$, and
    \item if $B\in\mathbb Q$, then $I(K_{\Ff}+B+\Mm_X)\sim 0$.
\end{enumerate}
\end{lem}
\begin{proof}
Since $B\not=0$ or $\Mm_X\not\equiv 0$, $K_{\Ff}$ is not pseudo-effective. By Theorem \ref{thm: subfoliation algebraic integrable}, there exists an algebraically integrable foliation $0\subsetneq\Ee\subset\Ff$. Since $\rk\Ff=1$, we deduce that $\rk\Ee=1$ and $\Ee=\Ff$. Thus $\Ff$ is algebraically integrable. By \cite[Theorem 3.10]{ACSS21}, possibly replacing $\Ii$ with $\Ii\cup\{1\}$ and $(X,\Ff,B)$ with a dlt model of $(X,\Ff,B)$, we may assume that $X$ is klt, $(X,\Ff,B)$ is dlt, and there exists a contraction $\pi: X\rightarrow Z$ such that $\Ff$ is induced by $\pi$. In particular, $\dim Z=1$. Let $F$ be a general fiber of $\pi$, $B_F:=B|_F$, and $\Mm^F:=\overline{\Mm_X|_F}$, then
$$K_F+B_F+\Mm^F_F=(K_X+B+\Mm_X)|_F=(K_{\Ff}+B+\Mm_X)|_F.$$
Since $\dim F=1,B\in\Ii$, and $p\Mm$ is base-point-free, there exists a finite set $\Ii_0\subset\Ii$ such that $B_F\in\Ii_0$. Since $(X,\Ff,B,\Mm)$ is lc, all components of $B$ are horizontal$/Z$. Thus $B\in\Ii_0$, and we get (1). 

Now we prove (2). By (1), we may assume that $\Ii=\Ii_0$ and $\Ii\subset\mathbb Q$. Since $\dim F=1,B\in\Ii$, and $p\Mm$ is base-point-free, there exists a positive integer $I$ depending only on $\Ii$ and $p$, such that $I(K_F+B_F+\Mm^F_F)\sim 0$. Thus $I(K_X+B+\Mm_X)\sim 0$ over the generic point of $Z$. We let
$$R:=\sum_D(\pi^*D-\pi^{-1}(D))$$
where the sum runs through all points on $Z$.

Since $Z$ is a curve and $K_{\Ff}+B+\Mm_X\sim_{\mathbb Q}0$, by Proposition \ref{prop: cbf is really cbf}(2.a), there exists a projective lc gfq $(Z,\Ff_Z,B_Z,\Mm^Z)$ induced by a canonical bundle formula of $\pi: (X,\Ff,B)\rightarrow Z$, such that $(Z,B_Z,\Mm^Z)$ is a g-pair induced by the canonical bundle formula of $\pi: (X,B-R,\Mm)\rightarrow Z$. Since $K_{\Ff}+B+\Mm_X\sim_{\mathbb Q}0$, $K_{\Ff_Z}+B_Z+\Mm^Z_Z\sim_{\mathbb Q}0$. Since $\Ff_Z$ is a foliation by points, $K_{\Ff_Z}=0$, $B_Z=0$, and $\Mm^Z_Z$ is torsion. Let $T:=\sum_{D\mid \pi^*D\not=\pi^{-1}(D)}\pi^{-1}(D)$ and $T_Z:=\sum_{D\mid \pi^*D\not=\pi^{-1}(D)}D$. By the definition of $R$ and $B_Z$,
$\left(X,B+T,\Mm\right)$
is lc and $(Z,T_Z,\Mm^Z)$ is a g-pair induced by the canonical bundle formula of $\pi: (X,B+T,\Mm)\rightarrow Z$.

If $\lfloor B\rfloor=0$, then by \cite[Chapter 6, Corollary 1]{Fil19}, possibly replacing $I$, we can choose $\Mm^Z$ so that $I\Mm^Z_Z$ is base-point-free. If $\lfloor B\rfloor\not=0$, then we let $S$ be an irreducible component of $\lfloor B\rfloor$, $(S,B_S,\Mm^S)$ the lc g-pair given by the adjunction
$$K_S+B_S+\Mm^S_S=(K_X+B+T+\Mm_X)|_S,$$
and $\pi_S: S\rightarrow Z$ the induced morphism. Then $p\Mm^S$ is base-point-free and $\pi_S$ is a finite morphism. Since $K_X+B+T+\Mm_X\sim_{\mathbb Q,Z}0$, by \cite[Theorem 1.1]{FS23}, $S|_F$ has at most two components, so $\deg(\pi_S)\leq 2$. Thus we may pick $\Mm^Z$ so that $2p\Mm^Z_Z$ is base-point-free. Possibly replacing $I$ with $2pI$, we may assume that $I\Mm^Z_Z$ is base-point-free. By Proposition \ref{prop: cbf m is generic},
$$I(K_X+B+\Mm_X)\sim I\pi^*(K_{\Ff_Z}+B_Z+\Mm^Z_Z)=I\pi^*\Mm^Z_Z=\pi^*(I\Mm^Z_Z)\sim 0.$$
\end{proof}

\begin{proof}[Proof of Theorem \ref{thm: surface index theorem global intro}]
If $\rk\Ff=0$, then $B=0$ and there is nothing left to prove. When $\rk\Ff=2$, the theorem follows from \cite[Proposition 7.2]{CH21}. So we may assume that $\rk\Ff=1$. Possibly replacing $\Ii$ with $\Ii\cup\{1\}$, we may assume that $1\in\Ii$. By Theorem \ref{thm: fdlt modification}, we may replace $(X,\Ff,B)$ with a dlt model and assume that $X$ is klt, $(X,\Ff,B)$ is dlt, and $\Ff$ is non-dicritical. In particular, $\Ff$ is canonical. If $B=0$, then the theorem follows from Lemma \ref{lem: surface canonical 12 index}. If $B\not=0$, then the theorem follows from Lemma \ref{lem: surface gloabl index with boundary}.
\end{proof}

The following result is an application of Theorem \ref{thm: surface index theorem global intro}. We do not need it in the rest of this paper.

\begin{thm}[Effective non-vanishing]
Let $\Ii\subset [0,1]\cap\mathbb Q$ be a DCC set. Then there exists a positive integer $I$ depending only on $\Ii$ satisfying the following. If $(X,\Ff,B)$ is a projective lc f-triple such that $\dim X=2$, $\rk\Ff=1$, $B\in\Ii$, and $\kappa_{\sigma}(K_{\Ff}+B)=0$, then $|I(K_{\Ff}+B)|\not=\emptyset$.
\end{thm}
\begin{proof}
By replacing $\Ii$ with $\Ii\cup\{1\}$, we may assume that $1\in\Ii$. By Theorem \ref{thm: fdlt modification}, we may replace $(X,\Ff,B)$ with a dlt model and assume that $X$ is $\mathbb Q$-factorial klt and $(X,\Ff,B)$ is dlt. By Theorems \ref{thm: can run mmp} and \ref{thm: abundance num dim 0}, we may run a $(K_{\Ff}+B)$-MMP which terminates with a good minimal model $(X',\Ff',B')$ of $(X,\Ff,B)$. Then $\kappa_{\sigma}(K_{\Ff'}+B')=\kappa(K_{\Ff'}+B')=0$. By \cite[Chapter V, Proposition 2.7(6)]{Nak04}, $K_{\Ff'}+B'\equiv 0$, so $K_{\Ff'}+B'\sim_{\mathbb Q}0$. By Theorem \ref{thm: surface index theorem global intro}, there exists a positive integer $I$ depending only on $\Ii$ such that $I(K_{\Ff'}+B')\sim 0$, hence $|I(K_{\Ff}+B)|\not=\emptyset$.
\end{proof}

\section{Proof of Theorem \ref{thm: global acc threefold}}\label{sec: Proof of Theorem}

\begin{prop}\label{prop: global acc with contraction}
Let $\Ii\subset [0,1]\cap\mathbb Q$ be a DCC set. Then there exists a positive integer $I$ depending only on $\Ii$ satisfying the following. Assume that
\begin{enumerate}
    \item $(X,\Ff,B)$ is a projective lc f-triple of dimension $\leq 3$,
    \item $B\in\Ii$,
    \item $\pi: X\rightarrow Z$ is a contraction such that $0<\dim Z<\dim X$,
    \item  there exists a foliation $\Ff_Z$ on $Z$ such that $\Ff=\pi^{-1}\Ff_Z$, and
    \item $K_{\Ff}+B\sim_{\mathbb Q}0$.
\end{enumerate} 
Then $I(K_{\Ff}+B)\sim 0$. In particular, the coefficients of $B$ belong to the finite set $\frac{1}{I}\mathbb N\cap [0,1]$.
\end{prop}
\begin{proof}
If $\rk\Ff=0$ then $B=0$, so we may take $I=1$ and we are done. If $\rk\Ff=\dim X$, then the proposition follows from \cite[Theorem 1.5]{HMX14} and \cite[Theorem 13]{Xu19}. Thus we may assume that $1\leq\rk\Ff\leq\dim X-1$. 

By Theorem \ref{thm: fdlt modification}, possibly replacing $\Ii$ with $\Ii\cup\{1\}$, we may assume that $(X,\Ff, B)$ is $\Qq$-factorial dlt. Let $F$ be a general fiber of $\pi$ and $B_F:=B|_F$. Then
 $$K_F+B_F=(K_X+B)|_F=(K_{\Ff}+B)|_F\sim_{\mathbb Q}0,$$
 and $B_F\in\Ii$. Since $(X,\Ff,B)$ is lc, $(F,B_F)$ is lc. By \cite[Theorem 1.5]{HMX14}, there exists a finite set $\Ii_0\subset\Ii$ depending only on $\Ii$ such that $B_F\in\Ii_0$.
 
By Proposition \ref{prop: cbf is really cbf}(2), there exists a projective lc gfq $(Z,\Ff_Z,B_Z,\Mm)$ induced by a canonical bundle formula of $\pi: (X,\Ff,B)\rightarrow Z$. Moreover, if $\dim X=3$ and $\dim Z=1$, then $\Ff_Z$ is a foliation by points, so $K_{\Ff_Z}=0$, $B_Z=0$, and $\Mm_Z\sim_{\mathbb Q}0$. Thus 
\begin{itemize}
    \item either $\dim Z=\dim X-1$, or
    \item $\dim X=3,\dim Z=1$, and $\Mm_Z\equiv 0$.
\end{itemize}
Since $K_{\Ff_Z}+B_Z+\Mm_Z\sim_{\mathbb Q}0$, if $\Ff_Z$ is a foliation by points, then $K_{\Ff_Z}=0$ and $B_Z=0$, an the proposition immediately follows from Proposition \ref{prop: low-dimension cbf}. Thus we may assume that $\Ff_Z$ is not a foliation by points, so $\dim X=3$ and $\dim Z=2$. By Proposition \ref{prop: low-dimension cbf} and Theorem \ref{thm: surface index theorem global intro}, there exists a positive integer $I_0$ depending only on $\Ii$, such that after possibly replacing $\Mm$, we have that $I_0(K_F+B_F)\sim 0$, $I_0\Mm$ is base-point-free, and $I_0(K_\Ff+B)\sim I_0\pi^*(K_{\Ff_Z}+B_Z+\Mm_Z)$. By \cite[Theorem 0.5]{Che22}, the coefficients of $B_Z$ belong to a DCC set depending only on $\Ii$. By Lemmas \ref{lem: surface canonical 12 index} and \ref{lem: surface gloabl index with boundary}, there exists a positive integer $I_1$ depending only on $\Ii$ such that $I_1(K_{\Ff_Z}+B_Z+\Mm_Z)\sim 0$. Thus
$$I_0I_1(K_{\Ff}+B)\sim I_0I_1\pi^*(K_{\Ff_Z}+B_Z+\Mm_Z)\sim 0,$$
so we may take $I=I_0I_1$.
\end{proof}

The proof of the following lemma is similar to \cite[Proof of Proposition 12.4]{CS21}.

\begin{lem}\label{lem: rk2 nonintegrable fibration}
Let $(X,\Ff,B)$ be a projective $\Qq$-factorial dlt f-triple such that $\dim X=3$ and $\rk\Ff=2$. Suppose that $K_{\Ff}+B\sim_{\mathbb R}0$, $B\not=0$, and $\Ff$ is not algebraically integrable. Let $B_0$ be an irreducible component of $B$. Then there exist a birational map $f: X\dashrightarrow X'$ and a contraction $\pi': X'\rightarrow Z$, such that 
\begin{enumerate}
    \item $f$ does not extract any divisor, 
    \item $B_0$ is not contracted by $f$,
    \item there exists a foliation $\Ff_Z$ on $Z$ such that $\Ff':=f_*\Ff=\pi'^{-1}\Ff_Z$, and
    \item $\dim Z=2$ and $Z$ is klt.
\end{enumerate}
\end{lem}
\begin{proof}
We run a $(K_{\Ff}+B-B\wedge B_0)$-MMP $f: X\dashrightarrow X'$ with scaling of an ample divisor $H$. By Theorem \ref{thm: can run mmp}, this MMP terminates with a Mori fiber space $\pi': X'\rightarrow Z$.  Let $\Ff':=f_*\Ff$ and let $B',B_0'$ be the images of $B,B_0$ on $X'$ respectively, and let $R$ be the $(K_{\Ff'}+B'-B'\wedge B_0')$-negative extremal ray contracted by $\pi': X'\rightarrow Z$. Then $\loc R=X'$.

(1) holds automatically. Since $K_{\Ff}+B\sim_{\mathbb R}0$, this MMP is also a $(-B_0)$-MMP, so $B_0$ is not contracted by $f$. Thus (2) holds. Since $\Ff'$ is not algebraically integrable, by \cite[Lemma 8.12 and Theorem 8.13]{Spi20}, $\dim Z=2$, and there exists a foliation $\Ff_Z$ on $Z$ such that $\Ff'=\pi'^{-1}\Ff_Z$, which implies (3).  Since $K_{X'}=K_{\Ff'}$ over the generic point of $Z$,
$$K_{X'}\cdot R=K_{\Ff'}\cdot R\leq (K_{\Ff'}+B'-B'\wedge B_0')\cdot R<0,$$
so $\pi'$ is a $K_{X'}$-Mori fiber space. Since $X'$ is klt, $Z$ is klt. This implies (4).
\end{proof}

\begin{proof}[Proof of Theorem \ref{thm: global acc threefold}]
By Theorem \ref{thm: surface index theorem global intro}, we may assume that $\dim X=3$. By \cite[Theorem 1.5]{HMX14}, \cite[Corollary 1.7]{Jia21}, and \cite[Theorem 13]{Xu19}, we may assume that $0<\rk\Ff<\dim X$. By Theorem \ref{thm: fdlt modification}, possibly replacing $\Ii$ with $\Ii\cup\{1\}$ and $(X,\Ff,B)$ with a dlt model, we may assume that $(X,\Ff,B)$ is $\Qq$-factorial dlt. If $B=0$, then there is nothing left to prove, so we may assume that $B\not=0$, hence $K_{\Ff}$ is not pseudo-effective. By Theorem \ref{thm: subfoliation algebraic integrable}, there exists an algebraically integrable foliation $0\not=\Ee\subset\Ff$. If $\Ff$ is algebraically integrable, then by \cite[Theorem 3.10]{ACSS21}, possibly replacing $(X,\Ff,B)$ with a dlt model, we may assume that there exists a contraction $\pi: X\rightarrow Z$ such that $\Ff$ is induced by $\pi$, and the theorem follows from Proposition \ref{prop: global acc with contraction}. Thus we may assume that $\Ff$ is not algebraically integrable. 

Let $S$ be a component of $B$. By Lemma \ref{lem: rk2 nonintegrable fibration}, there exists a birational map $f: X\dashrightarrow X'$ and a contraction $\pi': X'\rightarrow Z$ such that $f$ does not extract any divisor, $S$ is not contracted by $f$, $\Ff':=f_*\Ff$ is induced by a foliation $\Ff_Z$ on $Z$ (i.e. $\Ff'=\pi'^{-1}\Ff_Z$), and $\dim Z=2$. Let $B':=f_*B$ and $S':=f_*S$. Then $S'\not=0$. Since $K_{\Ff}+B\sim_{\mathbb Q}0$, $K_{\Ff'}+B'\sim_{\mathbb Q}0$ and $(X',\Ff',B')$ is lc. By Proposition \ref{prop: global acc with contraction}, there exists a positive integer $I$ depending only on $\Ii$ such that $I(K_{\Ff'}+B')\sim 0$. Thus $I(K_{\Ff}+B)\sim 0$, and the theorem follows.
\end{proof}

\end{document}